\documentclass[sn-mathphys,Numbered]{sn-jnl}

\usepackage{mathtools}
\usepackage{float}
\usepackage{helvet}
\usepackage{dsfont}
\usepackage[utf8]{inputenc}
\usepackage{makecell}
\usepackage{hyperref}
\usepackage{graphicx}
\usepackage{tikz}
\usepackage{pgfplots}
\pgfplotsset{compat=1.18}
\usetikzlibrary{decorations.pathmorphing,patterns}
\usetikzlibrary{shapes.multipart}
\usepackage{rotating}
\usepackage[font=footnotesize]{caption}
\usepackage{booktabs, makecell, caption}

\usepackage{graphicx}%
\usepackage{multirow}%
\usepackage{amsmath,amssymb,amsfonts}%
\usepackage{amsthm}%
\usepackage{mathrsfs}%
\usepackage[title]{appendix}%
\usepackage{xcolor}%
\usepackage{textcomp}%
\usepackage{manyfoot}%
\usepackage{booktabs}%
\usepackage{algorithm}%
\usepackage{algorithmicx}%
\usepackage{algpseudocode}%
\usepackage{listings}%
\usepackage{environ}
\usepackage{orcidlink}

\newcommand{\y}{\ensuremath{\mathbf{y}}}
\newcommand{\q}{\ensuremath{\mathbf{q}}}
\newcommand{\p}{\ensuremath{\mathbf{p}}}
\newcommand{\f}{\ensuremath{\mathbf{f}}}
\newcommand{\J}{\ensuremath{\mathbf{J}}}
\newcommand{\I}{\ensuremath{\mathbf{I}}}
\newcommand{\0}{\ensuremath{\mathbf{0}}}
\newcommand{\diag}{\ensuremath{\mathrm{diag}}}
\newcommand{\slow}{\ensuremath{\mathfrak{s}}}
\newcommand{\fast}{\ensuremath{\mathfrak{f}}}
\renewcommand{\H}{\ensuremath{\mathcal{H}}}

\makeatletter
\newsavebox{\measure@tikzpicture}
\NewEnviron{scaletikzpicturetowidth}[1]{%
  \def\tikz@width{#1}%
  \def\tikzscale{1}\begin{lrbox}{\measure@tikzpicture}%
  \BODY
  \end{lrbox}%
  \pgfmathparse{#1/\wd\measure@tikzpicture}%
  \edef\tikzscale{\pgfmathresult}%
  \BODY
}
\makeatother

\newtheorem{theorem}{Theorem}%
\newtheorem{example}{Example}%
\newtheorem{remark}{Remark}%
\newtheorem{lemma}{Lemma}

\newtheorem{definition}{Definition}%

\begin{document}

\title[Symplectic multirate generalized additive Runge--Kutta methods for Hamiltonian systems]{Symplectic multirate generalized additive Runge--Kutta methods for Hamiltonian systems}


\author*[1]{\fnm{Kevin} \sur{Schäfers} \orcidlink{0000-0002-6502-8106}}\email{schaefers@math.uni-wuppertal.de}

\author[1]{\fnm{Michael} \sur{Günther} \orcidlink{0000-0002-2195-4300}}\email{guenther@math.uni-wuppertal.de}

\author[2]{\fnm{Adrian} \sur{Sandu} \orcidlink{0000-0002-5380-0103}}\email{sandu@cs.vt.edu}

\affil[1]{\orgdiv{Institute of Mathematical Modelling, Analysis and Computational Mathematics (IMACM)}, \orgname{Bergische Universität Wuppertal}, \orgaddress{\street{Gau{\ss}strasse 20}, \city{Wuppertal}, \postcode{42119}, \country{Germany}}}

\affil[2]{\orgdiv{Computer Science Laboratory, Department of Computer Science}, \orgname{Virginia Polytechnic Institute and State University}, \orgaddress{\street{2202 Kraft Drive}, \city{Blacksburg}, \postcode{VA 24060},  \country{USA}}}


\abstract{The generalized additive Runge--Kutta (GARK) framework provides a powerful approach for solving additively partitioned ordinary differential equations. This work combines the ideas of symplectic GARK schemes and multirate GARK schemes to efficiently solve additively partitioned Hamiltonian systems with multiple time scales. Order conditions, as well as conditions for symplecticity and time-reversibility, are derived in the general setting of non-separable Hamiltonian systems. Investigations of the special case of separable Hamiltonian systems are also carried out. We show that particular partitions may introduce stability issues, and discuss partitions that enable an implicit-explicit integration leading to improved stability properties. Higher-order symplectic multirate GARK schemes based on advanced composition techniques are discussed. The performance of the schemes is demonstrated by means of the Fermi--Pasta--Ulam problem.}




\maketitle

\section{Introduction}
The numerical simulation of Hamiltonian systems with multiscale dynamics is challenging. Small step sizes $h$ are required to capture the system's fast dynamics. On the other hand, the simulation of the slow dynamics with a larger step size $H=h\cdot M, \ M \in \mathbb{N}, \ M \gg 1$, is sufficiently accurate. In many applications, e.g. in lattice quantum chromodynamics \cite{knechtli2017lattice}, the slow subsystem is characterized by expensive evaluation costs. Hence the use of single rate integration schemes using the small step size $h$ for the entire Hamiltonian system results in an inefficient integration because of unnecessary costly function evaluations. 
Another example of Hamiltonian systems showing dynamics on different time scales is given by mechanical systems with stiff and non-stiff components. This work aims at deriving multirate integrators for efficient and structure-preserving simulation of these Hamiltonian systems. 

\subsection{Geometric Integration}

Consider a twice continuously differentiable Hamiltonian function $\H \, : \, \mathbb{R}^{n_p} \times \mathbb{R}^{n_q} \to \mathbb{R}$ that splits into $N\in \mathbb{N}$ partitions
\begin{align}\label{eq:Hamiltonian}
    \H(\p,\q) = \sum\limits_{m=1}^N \H^{\{m\}}(\p,\q),
\end{align}
 where $\q \in \mathbb{R}^{n_q}$ denote the generalized coordinates and $\p \in \mathbb{R}^{n_p}$ the conjugate momenta ($n_p = n_q = n/2$).
Based on this Hamiltonian function, we are able to formulate the Hamiltonian equations of motion as an additively partitioned initial value problem (IVP) of ordinary differential equations (ODEs)
\begin{equation}\label{eq:EoM}
    \dot{\y} = \J^{-1} \nabla \H(\y) = \J^{-1} \sum\limits_{m=1}^N \nabla \H^{\{m\}}(\y), \quad \y(t_0)= \y_0,
\end{equation}
with $\y = (\p,\q)^\top$ and 
\begin{equation}\label{eq:structure_matrix}
    \J = \begin{pmatrix}
\0_{n_p \times n_p} & \I_{n_p \times n_q} \\
-\I_{n_q \times n_p} & \0_{n_q \times n_q}
\end{pmatrix}.
\end{equation}
The Hamiltonian flow $\y(t) = \varphi_t(\y_0)$, i.e., the exact solution to \eqref{eq:EoM}, is characterized by the following properties \cite{HairerLubichWanner}:

\begin{itemize}

    \item \textbf{Energy conservation}: The Hamiltonian is an invariant of the flow:
    \begin{align}\label{eq:energy_conservation}
        \frac{d}{dt} \H(\varphi_t(\y_0)) = 0.
    \end{align}
    
    \item \textbf{$\rho$-reversibility}: The Hamiltonian is invariant with respect to changing the sign of momenta, $\H(\p,\q) = \H(-\p,\q)$, which can be formalized as \vspace*{-1ex}
    \begin{align}\label{eq:rho_definition}
        \H = \H \circ \rho, \quad \rho \coloneqq \begin{pmatrix}
            -\mathbf{I}_{n_p \times n_p} & \mathbf{0}_{n_p \times n_q} \\
            \hphantom{-}\mathbf{0}_{n_q \times n_p} & \mathbf{I}_{n_q \times n_q}.
        \end{pmatrix}
    \end{align}
    Consequently, the Hamiltonian equations of motion \eqref{eq:EoM} are $\rho$-reversible, i.e.,     
    \begin{align}\label{eq:rho-reversibility}
        \rho \circ (\nabla \H) = - \nabla (\H \circ \rho).
    \end{align}
    
    \item \textbf{Time-reversibility}: The Hamiltonian flow is time-reversible, i.e., \begin{align}\label{eq:time-reversibility}
        \rho \circ \varphi_t \circ \rho \circ \varphi_t(\y_0) = \y_0 \quad \Leftrightarrow \quad \rho \circ \varphi_t = \varphi_{-t} \circ \rho,
    \end{align}
    where the second equivalent equation hold due to the symmetry of the flow
    \begin{align}\label{eq:symmetry}
        \varphi_t \circ \varphi_{-t}(\y_0) = \y_0.
    \end{align}
    \item\textbf{Symplecticity}: The Hamiltonian flow is symplectic, i.e., \begin{align}\label{eq:symplecticity}
        \left(\frac{\partial \varphi_t (\y_0)}{\partial \y_0}\right)^\top \J \left( \frac{\partial \varphi_t(\y_0)}{\partial \y_0}\right) = \J.
    \end{align}
\end{itemize}

We demand the numerical integration scheme $\Phi_H : \mathbb{R}^n \to \mathbb{R}^n, \  \Phi_H(\y_0) = \y_1 \approx \varphi_H(\y_0)$ to preserve the time-reversibility \eqref{eq:time-reversibility}, as well as the symplecticity \eqref{eq:symplecticity} of the Hamiltonian flow. For time-reversibility, we obtain the criterion
\begin{align}\label{eq:discrete_time-reversibility}
    \rho \circ \Phi_H \circ \rho \circ \Phi_H(\y_0) = \y_0, 
\end{align}
whereas for symplecticity one obtains
\begin{align}\label{eq:discrete_symplecticity}
    \left(\frac{\partial \Phi_H(\y_0)}{\partial \y_0}\right)^\top \J \left(\frac{\partial \Phi_H(\y_0)}{\partial \y_0}\right) =  \J.
\end{align}

\noindent Performing a backward error analysis, it can be shown that symplectic schemes preserve a shadow Hamiltonian $\tilde{\H}$ that is close to $\H$ \cite{HairerLubichWanner}. Consequently, the total energy oscillates close to the real value, explaining their good long-term energy behaviour. In contrast, energy-conserving integration schemes preserve the Hamiltonian up to numerical accuracy \cite{HairerLubichWanner}. It is well known that numerical time integration schemes using constant time steps cannot be both symplectic and energy-conserving \cite{chartier2006algebraic}.

There are many approaches for symplectic integration of Hamiltonian systems like decomposition methods \cite{omelyan2003symplectic}, symplectic Runge--Kutta methods \cite{sanz2018numerical} and variational integrators \cite{marsden2001discrete}. In this work, we focus on symplectic generalized additive Runge--Kutta (GARK) methods \cite{gunther2021symplectic} which are based on a generalized-structure approach to additive Runge--Kutta methods \cite{sandu2015generalized}.

\subsection{Multirate Integration}

Multirate integration for ODE systems exhibiting multiscale dynamics uses different time steps to solve different time scales, e.g. \cite{gunther2001multirate,gear1984multirate,savcenco2007multirate,bartel2020inter,schafers2023spline}. However, these multirate approaches do not focus on the preservation of geometric properties of the underlying system. For the special case of symplectic systems where the slow potential is costly to evaluate and the fast potential is cheap to evaluate, multirate integrators within the framework of splitting and composition methods have been successful. For example, the idea of impulse methods \cite{HairerLubichWanner} belongs to this framework. Moreover, multirate integration based on splitting and composition methods has become very popular in molecular dynamics and in particular finds application in lattice quantum chromodynamics \cite{clark2011improving} (in the physics literature, it is referred to nested integration). 
Symmetric compositions of symplectic schemes provide a tool for multirate integration of Hamiltonian systems of arbitrarily high convergence order. Alternatively, an approach based on the variational principle results in multirate integration schemes for constrained systems \cite{leyendecker2013variational}. This work focuses on the derivation of structure-preserving multirate schemes in the GARK framework. Multirate GARK schemes have already been introduced in \cite{gunther2016multirate}.

\subsection{Stability}

Symplectic integrators are not energy-preserving, but they preserve exactly a shadow Hamiltonian \cite{skeel2001practical,kennedy2013shadow} which is close to the original Hamiltonian, explaining their excellent long-time behavior. For stiff Hamiltonian problems, one also has to pay attention to the stability properties of the numerical integration scheme \cite{mclachlan2004nonlinear}. The GARK framework provides a theory for linear and nonlinear stability \cite{sandu2015generalized}.

This paper develops the new subclass of symplectic and symmetric multirate GARK (sMGARK) schemes, combining the ideas of symplectic GARK schemes and multirate GARK schemes. It provides a framework for geometric multirate integration of Hamiltonian systems that includes an order theory, as well as conditions for nonlinear stability \cite{sandu2015generalized}.
The remaining part of the paper is organized as follows.
Section \ref{sec:sMGARK} derives the subclass of sMGARK schemes in the general setting of non-separable Hamiltonian splittings. Order conditions up to order three, as well as conditions for symplecticity and time-reversibility are derived. Furthermore, we show that all schemes in this subclass are nonlinearly stable. 
Section \ref{sec:P-sMGARK} considers the important special case of separable Hamiltonian splittings. Here, we introduce partitioned sMGARK methods which exploit the special structure of the partitions and allow for the construction of explicit symplectic schemes. The explicitness of the integrator introduces stability concerns which one may overcome by using different partitions that enable the derivation of implicit-explicit integration schemes. 
These stability-enhanced methods are discussed in Section \ref{sec:IMEX-sMGARK}. 
Section \ref{sec:composition} is devoted to the derivation of efficient (partitioned) sMGARK schemes of higher order based on (advanced) composition techniques \cite{suzuki1990fractal,yoshida1990construction,omelyan2002construction,kahan1997composition}. 
In Section \ref{sec:numerics}, the performance of the schemes is demonstrated by means of numerical tests on the Fermi--Pasta--Ulam (FPU) problem, a standard benchmark problem for multirate integration of Hamiltonian systems.
Section \ref{sec:conclusions} draws conclusions and points to future research directions.

\section{sMGARK schemes for general splittings}\label{sec:sMGARK}
Consider a general (possibly non-separable) Hamiltonian splitting 
\begin{align}\label{eqref:Hamiltonian_splitting}
    \mathcal{H}(\mathbf{y}) = \sum\limits_{m=1}^N \mathcal{H}^{\{m\}}(\mathbf{y})
\end{align}
and the corresponding equations of motion 
\begin{align}\label{eq:IVP}
    \dot{\mathbf{y}} = \mathbf{J}^{-1} \left( \sum\limits_{m=1}^N \nabla \mathcal{H}^{\{m\}}(\mathbf{y})\right) \eqqcolon \sum\limits_{m=1}^N \mathbf{f}^{\{m\}}(\mathbf{y}), \quad \mathbf{y}(t_0) = \mathbf{y}_0, 
\end{align}
that define an additively partitioned ODE system. $\mathbf{J}$ is assumed to be a skew-symmetric, regular matrix but does not necessarily need to coincide with \eqref{eq:structure_matrix}.
We start with a brief recapitulation of generalized additive Runge--Kutta (GARK) methods \cite{sandu2015generalized}. One step of a GARK scheme, applied to the additively partitioned system \eqref{eq:IVP} advances the solution $\mathbf{y}_0$ at $t_0$ to the solution $\mathbf{y}_1$ at $t_1 = t_0 + H$ as follows:
\begin{subequations}\label{eq:GARK}
\begin{align}
    \y_{1} &= \y_0 + H \sum\limits_{q=1}^N \sum\limits_{i=1}^{s^{\{q\}}} b_i^{\{q\}} \f^{\{q\}}\left( \mathbf{Y}_i^{\{q\}}\right), \\
    \mathbf{Y}_i^{\{q\}} &= \y_0 + H \sum\limits_{m=1}^N \sum\limits_{j=1}^{s^{\{m\}}} a_{i,j}^{\{q,m\}} \f^{\{m\}}\left(\mathbf{Y}_j^{\{m\}}\right), \quad i=1,\ldots,s^{\{q\}}, \, q=1,\ldots,N.  
\end{align}
\end{subequations}
The corresponding generalized Butcher tableau is
\begin{align}\label{eq:GARK_Butcher}
    \begin{tabular}{c c c c}
        $\mathbf{A}^{\{1,1\}}$ & $\mathbf{A}^{\{1,2\}}$ & $\cdots$ & $\mathbf{A}^{\{1,N\}}$ \\
        $\mathbf{A}^{\{2,1\}}$ & $\mathbf{A}^{\{2,2\}}$ & $\cdots$ & $\mathbf{A}^{\{2,N\}}$ \\
        $\vdots$ & $\vdots$ & & $\vdots$ \\
        $\mathbf{A}^{\{N,1\}}$ & $\mathbf{A}^{\{N,2\}}$ & $\cdots$ & $\mathbf{A}^{\{N,N\}}$ \\\Xhline{2\arrayrulewidth}\rule{0ex}{3ex} 
        $\mathbf{b}^{\{1\}\top}$ & $\mathbf{b}^{\{2\}\top}$ & $\cdots$ & $\mathbf{b}^{\{N\}\top}$
    \end{tabular}\,.
\end{align}
Whereas additive Runge--Kutta (ARK) schemes \cite{kennedy2003additive} use the same stage value for different components of the right-hand side, GARK schemes allow for using different stage values for different components. The methods defined by the Butcher table $(\mathbf{A}^{\{q,q\}},\mathbf{b}^{\{q\}})$, $q=1,\ldots,N$, can be viewed as stand-alone integration schemes applied to each individual component $q$. The off-diagonal matrices $\mathbf{A}^{\{q,m\}}$, $q,m=1,\ldots,N$, $q \neq m$, define the coupling among the components $q$ and $m$.
Based on the coefficients in \eqref{eq:GARK_Butcher}, we define 
\begin{subequations}\label{eq:aux_matrices}
\begin{align}
    \mathbf{c}^{\{q,m\}} &\coloneqq \mathbf{A}^{\{q,m\}} \cdot \mathbf{1} \in \mathbb{R}^{s^{\{q\}}}, \quad \mathbf{1} \coloneqq \left( 1 \ \cdots \ 1\right)^\top \in \mathbb{R}^{s^{\{m\}}}, \\
    \mathbf{B}^{\{q\}} &\coloneqq \mathrm{diag}\left(\mathbf{b}^{\{q\}}\right) \in \mathbb{R}^{s^{\{q\}} \times s^{\{q\}}}, \\
    \mathbf{P}^{\{q,m\}} &\coloneqq \mathbf{A}^{\{m,q\}\top} \mathbf{B}^{\{m\}} + \mathbf{B}^{\{q\}} \mathbf{A}^{\{q,m\}} - \mathbf{b}^{\{q\}} \mathbf{b}^{\{m\}\top} \in \mathbb{R}^{s^{\{q\}} \times s^{\{m\}}}, \\ \label{eq:definition_P}
    \mathbf{P} &\coloneqq \left( \mathbf{P}^{\{q,m\}}\right)_{1 \leq q,m \leq N} \in \mathbb{R}^{s \times s}, \quad s \coloneqq \sum_{q=1}^N s^{\{q\}},
\end{align}
\end{subequations}
for all $q,m = 1,\ldots,N$. 
Multirate GARK methods \cite{gunther2016multirate} are a subclass of GARK schemes that are able to exploit multiscale behavior. Consider a two-way partitioned Hamiltonian system 
\begin{align}\label{eq:MR-EoM}
       \dot{\y} &= \J^{-1} \left(\nabla \H^{\{\slow\}}(\y) + \nabla \H^{\{\fast\}}(\y)\right) \eqqcolon \f^{\{\slow\}}(\y) + \f^{\{\fast\}}(\y), \quad \y(t_0) = \y_0,
\end{align}
with one slow component $\slow$ and one fast component $\fast$. Here, we use a large macro-step size $H$ for the slow component, and a small micro-step size $h= H/M$ for the fast part, where $M \in \mathbb{N}$ is called the \textit{multirate factor}.
One macro-step of a MGARK scheme, applied to the two-way partitioned system \eqref{eq:MR-EoM}, advances the solution $\y_0$ at $t_0$ to the solution $\y_1$ at time point $t_1 = t_0 + H$ with $M \in \mathbb{N}$ equal micro-steps of size $h=H/M$ as follows:
\allowdisplaybreaks
\begin{subequations}\label{eq:MGARK}
\begin{align}
    \y_1 = \y_0 &+ h \sum\limits_{\lambda = 1}^M \sum\limits_{i=1}^{\bar{s}^{\{\fast,\lambda\}}} \bar{b}_i^{\{\fast,\lambda\}} \f^{\{\fast\}}\left(\mathbf{Y}_i^{\{\fast,\lambda\}}\right) + H \sum\limits_{i=1}^{\bar{s}^{\{\slow\}}} \bar{b}_i^{\{\slow\}} \f^{\{\slow\}}\left(\mathbf{Y}_i^{\{\slow\}}\right), \\
    \begin{split}
    \mathbf{Y}_i^{\{\slow\}} = \y_0 &+ H\!\sum\limits_{j=1}^{\bar{s}^{\{\slow\}}} \bar{a}_{i,j}^{\{\slow,\slow\}} \f^{\{\slow\}}\!\left(\mathbf{Y}_j^{\{\slow\}}\right) \\
    &+ h \sum\limits_{\lambda = 1}^M \sum\limits_{j=1}^{\bar{s}^{\{\fast,\lambda\}}} \bar{a}_{i,j}^{\{\slow,\fast,\lambda\}} \f^{\{\fast\}}\!\left(\mathbf{Y}_j^{\{\fast,\lambda\}}\right), \ i=1,\ldots,\bar{s}^{\{\slow\}}, 
    \end{split} \\
    \begin{split}
    \mathbf{Y}_i^{\{\fast,\lambda\}} = \y_0 &+ h \sum\limits_{\ell=1}^{\lambda-1} \sum\limits_{j=1}^{\bar{s}^{\{\fast,\ell\}}} \bar{b}_j^{\{\fast,\ell\}} \f^{\{\fast\}}\left(\mathbf{Y}_j^{\{\fast,\ell\}}\right) + H \sum \limits_{j=1}^{\bar{s}^{\{\slow\}}} \bar{a}_{i,j}^{\{\fast,\slow,\lambda\}} \f^{\{\slow\}}\left(\mathbf{Y}_j^{\{\slow\}}\right) \\
    &+ h \sum\limits_{j=1}^{\bar{s}^{\{\fast,\lambda\}}} \bar{a}_{i,j}^{\{\fast,\fast,\lambda\}} \f^{\{\fast\}}\left( \mathbf{Y}_j^{\{\fast,\lambda\}}\right), \  i=1,\ldots,\bar{s}^{\{\fast,\lambda\}}, \quad \lambda = 1,\ldots,M. 
    \end{split}
\end{align}
\end{subequations}
The base schemes are given by Runge--Kutta methods $(\mathbf{\bar{A}}^{\{\slow,\slow\}},\mathbf{\bar{b}}^{\{\slow\}})$ for the slow subsystem and $(\mathbf{\bar{A}}^{\{\fast,\fast,\lambda\}},\mathbf{\bar{b}}^{\{\fast,\lambda\}}),\ \lambda\!=\!1,\ldots,M,$ for the fast subsystem. The coefficients $\mathbf{\bar{A}}^{\{\slow,\fast,\lambda\}}$, $\mathbf{\bar{A}}^{\{\fast,\slow,\lambda\}}$, $\lambda\!=\!1,\ldots,M$, realize the coupling of the $\lambda$-th micro-step with the macro-step. The MGARK scheme \eqref{eq:MGARK} can be written as a GARK scheme \eqref{eq:GARK} over the macro-step $H$ with fast stage vectors $\mathbf{Y}^{\{\fast\}} \coloneqq \left(\mathbf{Y}^{\{\fast,1\}\top}\  \cdots\  \mathbf{Y}^{\{\fast,M\}\top}\right)^\top$. Its generalized Butcher tableau \eqref{eq:GARK_Butcher} reads
\setlength{\tabcolsep}{0.4ex}
\begin{align}\label{eq:MGARK_Butcher}
    \raisebox{-25.25pt}{$\renewcommand{\arraystretch}{1.5} \begin{tabular}{c|c} 
         $\mathbf{A}^{\{\fast,\fast\}}$ & $\mathbf{A}^{\{\fast,\slow\}}$ \\\hline\rule{0pt}{3ex}
         $\mathbf{A}^{\{\slow,\fast\}}$ & $\mathbf{A}^{\{\slow,\slow\}}$ \\\Xhline{2\arrayrulewidth}
         $\mathbf{b}^{\{\fast\}\top}$ & $\mathbf{b}^{\{\slow\}\top}$
    \end{tabular} 
    \renewcommand{\arraystretch}{1} $}
    \!\raisebox{-34.75pt}{$~\coloneqq~$}\!
     \begin{tabular}{c c c c | c}
         $\frac{1}{M} \mathbf{\bar{A}}^{\{\fast,\fast,1\}}$ & $\mathbf{0}$ & $\cdots$ & $\mathbf{0}$ & $\mathbf{\bar{A}}^{\{\fast,\slow,1\}}$ \\
         $\frac{1}{M} \mathbf{1} \mathbf{\bar{b}}^{\{\fast,1\}\top}$ & $\frac{1}{M} \mathbf{\bar{A}}^{\{\fast,\fast,2\}}$ & $\ddots$ & $\vdots$ & $\mathbf{\bar{A}}^{\{\fast,\slow,2\}}$ \\
         $\vdots$ & $\ddots$ & $\ddots$ & $\mathbf{0}$ & $\vdots$ \\
         $\frac{1}{M} \mathbf{1} \mathbf{\bar{b}}^{\{\fast,1\}\top}$ & $\cdots$ & $\frac{1}{M} \mathbf{1} \mathbf{\bar{b}}^{\{\fast,M-1\}\top}$ & $\frac{1}{M} \mathbf{\bar{A}}^{\{\fast,\fast,M\}}$ & $\mathbf{\bar{A}}^{\{\fast,\slow,M\}}$ \\[-2ex]
         & & & & 
         \\\hline\rule{0pt}{3ex}
         $\frac{1}{M} \mathbf{\bar{A}}^{\{\slow,\fast,1\}}$ & $\frac{1}{M} \mathbf{\bar{A}}^{\{\slow,\fast,2\}}$ & $\cdots$ & $\frac{1}{M} \mathbf{\bar{A}}^{\{\slow,\fast,M\}}$ & $\mathbf{\bar{A}}^{\{\slow,\slow\}}$ \\[-2ex]
         & & & & \\ \Xhline{2\arrayrulewidth}\rule{0pt}{3ex}
         $\frac{1}{M} \mathbf{\bar{b}}^{\{\fast,1\}\top}$ & $\frac{1}{M} \mathbf{\bar{b}}^{\{\fast,2\}\top}$ & $\cdots$ & $\frac{1}{M} \mathbf{\bar{b}}^{\{\fast,M\}\top}$ & $\mathbf{\bar{b}}^{\{\slow\}\top}$
    \end{tabular}\,.
\end{align}
In \cite{gunther2021symplectic}, conditions for symmetry and symplecticity of GARK schemes have been derived. As MGARK schemes are just a subclass of GARK schemes, we can apply the conditions in a straightforward manner to the special structure of MGARK schemes. After removing redundant conditions, one obtains the following.\\

\definition[Symmetric MGARK schemes]{
Let $\mathcal{P}$ denote permutation matrices of matching dimensions that reverse the entries of a vector. The MGARK scheme \eqref{eq:MGARK}--\eqref{eq:MGARK_Butcher} is symmetric if for all $\lambda=1,\ldots,M$ it holds
\begin{align}
    \mathbf{\bar{b}}^{\{\slow\}} &= \mathcal{P} \mathbf{\bar{b}}^{\{\slow\}}, &
    \mathbf{\bar{A}}^{\{\slow,\slow\}} &= \mathbf{1} \mathbf{\bar{b}}^{\{\slow\}\top} \!-\! \mathcal{P} \mathbf{\bar{A}}^{\{\slow,\slow\}} \mathcal{P}, \nonumber \\
    \mathbf{\bar{b}}^{\{\fast,\lambda\}} &=\mathcal{P}\mathbf{\bar{b}}^{\{\fast,M+1-\lambda\}}, &
    \mathbf{\bar{A}}^{\{\fast,\fast,\lambda\}} &= \mathbf{1} \mathbf{\bar{b}}^{\{\fast,\lambda\}\top} \!-\! \mathcal{P}\mathbf{\bar{A}}^{\{\fast,\fast,M+1-\lambda\}} \mathcal{P}, \label{eq:symmetric_MGARK} \\
    \mathbf{\bar{A}}^{\{\fast,\slow,\lambda\}} &= \mathbf{1} \mathbf{\bar{b}}^{\{\slow\}\top} \!-\! \mathcal{P} \mathbf{\bar{A}}^{\{\fast,\slow,M+1-\lambda\}} \mathcal{P}, &
    \mathbf{\bar{A}}^{\{\slow,\fast,\lambda\}} &= \mathbf{1} \mathbf{\bar{b}}^{\{\fast,\lambda\}\top} \!-\! \mathcal{P} \mathbf{\bar{A}}^{\{\slow,\fast,M+1-\lambda\}} \mathcal{P}. \nonumber
\end{align}
}\normalfont
\remark{As in \eqref{eqref:Hamiltonian_splitting} all partitions are assumed to be Hamiltonian, symmetric MGARK schemes are also time-reversible \cite{gunther2021symplectic}.}\\

\definition[Symplectic MGARK schemes]{
Consider a MGARK scheme \eqref{eq:MGARK} given by the generalized Butcher tableau \eqref{eq:MGARK_Butcher}. The method is symplectic if
\begin{subequations}\label{eq:symplectic_MGARK}
\begin{align}
    \mathbf{\bar{A}}^{\{\slow,\slow\}\top} \mathbf{\bar{B}}^{\{\slow\}} + \mathbf{\bar{B}}^{\{\slow\}} \mathbf{\bar{A}}^{\{\slow,\slow\}} - \mathbf{\bar{b}}^{\{\slow\}} \mathbf{\bar{b}}^{\{\slow\}\top} &= \mathbf{0}, \label{eq:symplectic_MGARK_a}
\intertext{and for all $\lambda=1,\ldots,M$ it holds}
    \mathbf{\bar{A}}^{\{\fast,\fast,\lambda\}\top} \mathbf{\bar{B}}^{\{\fast,\lambda\}} + \mathbf{\bar{B}}^{\{\fast,\lambda\}} \mathbf{\bar{A}}^{\{\fast,\fast,\lambda\}} - \mathbf{\bar{b}}^{\{\fast,\lambda\}} \mathbf{\bar{b}}^{\{\fast,\lambda\}\top} &= \mathbf{0}, \label{eq:symplectic_MGARK_b} \\ 
    \mathbf{\bar{A}}^{\{\fast,\slow,\lambda\}\top} \mathbf{\bar{B}}^{\{\fast,\lambda\}} + \mathbf{\bar{B}}^{\{\slow\}} \mathbf{\bar{A}}^{\{\slow,\fast,\lambda\}} - \mathbf{\bar{b}}^{\{\slow\}} \mathbf{\bar{b}}^{\{\fast,\lambda\}\top} &= \mathbf{0}. \label{eq:symplectic_MGARK_c}
\end{align}
\end{subequations}
}\normalfont 
The well-known property that symplectic Runge--Kutta schemes with positive weights $b_i > 0$ are algebraically stable \cite{burrage1979stability} generalizes to the MGARK setting, as the following theorem states.\\
\theorem[Algebraic stability of sMGARK schemes]{
Symplectic MGARK schemes are algebraically stable, provided that $\bar{b}_i^{\{\slow\}} > 0, \ i=1,\ldots,\bar{s}^{\{\slow\}},$ and $\bar{b}_i^{\{\fast,\lambda\}} > 0, \ i=1,\ldots,\bar{s}^{\{\fast,\lambda\}}$ for all $\lambda=1,\ldots,M.$
}
\proof{
Consider a symplectic MGARK method. Then, the coupling parts are symplectic, i.e., \eqref{eq:symplectic_MGARK_c} is satisfied. This is equivalent to 
$\mathbf{P}^{\{q,m\}} = \mathbf{0}, \ q \neq m$.
According to \cite{sandu2015generalized}, a GARK scheme with this property is stability-decoupled, i.e., the interaction of different components does not influence the overall nonlinear stability. Hence it remains to show that the Runge--Kutta base schemes $(\mathbf{\bar{A}}^{\{\slow,\slow\}},\mathbf{\bar{b}}^{\{\slow\}})$ and $(\mathbf{\bar{A}}^{\{\fast,\fast,\lambda\}},\mathbf{\bar{b}}^{\{\fast,\lambda\}}),\ \lambda=1,\ldots,M,$ are algebraically stable. As the MGARK schemes is symplectic, it holds \eqref{eq:symplectic_MGARK_a}--\eqref{eq:symplectic_MGARK_b} and thus the Runge--Kutta base schemes are symplectic. Moreover, all weights are assumed to be positive, resulting in algebraically stable base schemes \cite{burrage1979stability}. Hence the overall scheme is algebraically stable. \hfill \endproof
}

\remark[Nonlinear stability]{Following the proofs in \cite{sandu2015generalized}, an algebraically stable sMGARK scheme \eqref{eq:MGARK} applied to \eqref{eq:IVP} with dispersive component functions \begin{align*}
    \left\langle \f^{\{m\}}(\y) - \f^{\{m\}}(\mathbf{\tilde{y}}), \y - \mathbf{\tilde{y}} \right\rangle \leq \nu^{\{m\}} \lVert \y - \mathbf{\tilde{y}} \rVert^2, \quad \nu^{\{m\}} <0,
\end{align*} 
is unconditionally nonlinearly stable, in the sense that the difference of any two numerical solutions is non-increasing, i.e., it holds $\lVert \Delta \y_1\rVert \leq \lVert \Delta \y_0\rVert$ for all step sizes $H>0$.}\normalfont\\

The GARK framework provides an order theory based on N-trees \cite{araujo1997symplectic}. In \cite{gunther2016multirate}, the order conditions have been adapted to the MGARK setting. The order conditions up to order three, excluding redundant ones, are given in Tables \ref{tab:MGARK_order_slow} and \ref{tab:MGARK_order_fast}. Here we used the common notation of the identity matrix $\mathbf{I}$, as well as the vectors $\mathbf{1} = (1\ \cdots\ 1)^\top$.  

In case of symmetric and/or symplectic MGARK schemes, the number of order conditions reduces significantly. As symmetric schemes are of even order, the order conditions for even $p$ are automatically satisfied. In case of symplectic MGARK schemes \eqref{eq:symplectic_MGARK} we have the following result.\\

\begin{subequations}\label{eq:order_slow}
\noindent\begin{minipage}{\linewidth}
    \newcommand{\stepequation}{\refstepcounter{equation}\thetag{\theequation}}
    \captionof{table}{Slow order conditions for the MGARK scheme \eqref{eq:MGARK}-\eqref{eq:MGARK_Butcher}}
    \label{tab:MGARK_order_slow}
    \begin{tabular}{p{.05\textwidth} p{.9\textwidth}}
    \toprule
    $p$ & Order condition  \\
    \midrule 
    1 & $\mathbf{\bar{b}}^{\{\slow\}\top}\mathbf{1} = 1$ \hfill\stepequation\hspace*{.5ex}\label{eq:order_slow_1a}\\
    2 & $\mathbf{\bar{b}}^{\{\slow\}\top} \mathbf{\bar{A}}^{\{\slow,\slow\}}\mathbf{1} = \tfrac{1}{2} \hfill\stepequation\hspace*{.5ex}\label{eq:order_slow_2a}$ \\
     & $\mathbf{\bar{b}}^{\{\slow\}\top} \left(\sum_{\lambda=1}^M \mathbf{\bar{A}}^{\{\slow,\fast,\lambda\}} \mathbf{1}\right) = \frac{M}{2}\hfill{\stepequation\hspace*{.5ex}\label{eq:order_slow_2b}}$  \\
     3 & $\mathbf{\bar{b}}^{\{\slow\}\top} \diag\left(\mathbf{\bar{A}}^{\{\slow,\slow\}}\mathbf{1}\right)\mathbf{\bar{A}}^{\{\slow,\slow\}}\mathbf{1} = \frac{1}{3}\hfill\stepequation \hspace*{.5ex}\label{eq:order_slow_3a}$ \\
    & $\mathbf{\bar{b}}^{\{\slow\}\top} \diag\left(\mathbf{\bar{A}}^{\{\slow,\slow\}}\mathbf{1}\right) \left(\sum_{\lambda=1}^M \mathbf{\bar{A}}^{\{\slow,\fast,\lambda\}}\mathbf{1}\right) = \frac{M}{3}$ \hfill\refstepcounter{equation}\thetag{\theequation} \label{eq:order_slow_3b}\\
     & $\mathbf{\bar{b}}^{\{\slow\}\top} \diag\left(\sum_{\lambda=1}^M \mathbf{\bar{A}}^{\{\slow,\fast,\lambda\}}\mathbf{1}\right) \left(\sum_{\lambda=1}^M \mathbf{\bar{A}}^{\{\slow,\fast,\lambda\}}\mathbf{1}\right) = \frac{M^2}{3}$ \hfill\refstepcounter{equation}\thetag{\theequation} \label{eq:order_slow_3c}\\
     & $\mathbf{\bar{b}}^{\{\slow\}\top} \mathbf{\bar{A}}^{\{\slow,\slow\}} \mathbf{\bar{A}}^{\{\slow,\slow\}} \mathbf{1} = \frac{1}{6}$ \hfill\refstepcounter{equation}\thetag{\theequation} \label{eq:order_slow_3d}\\ 
     & $\mathbf{\bar{b}}^{\{\slow\}\top} \mathbf{\bar{A}}^{\{\slow,\slow\}} \left(\sum_{\lambda=1}^M \mathbf{\bar{A}}^{\{\slow,\fast,\lambda\}}\mathbf{1}\right) = \frac{M}{6}$ \hfill\refstepcounter{equation}\thetag{\theequation} \label{eq:order_slow_3e} \\
     & $\mathbf{\bar{b}}^{\{\slow\}\top} \left(\sum_{\lambda=1}^M  \mathbf{\bar{A}}^{\{\slow,\fast,\lambda\}} \mathbf{\bar{A}}^{\{\fast,\slow,\lambda\}}\right)\mathbf{1} = \frac{M}{6}$ \hfill\refstepcounter{equation}\thetag{\theequation} \label{eq:order_slow_3f} \\
     & $\mathbf{\bar{b}}^{\{\slow\}\top} \left(\sum_{\lambda=1}^M \mathbf{\bar{A}}^{\{\slow,\fast,\lambda\}} ( \mathbf{\bar{A}}^{\{\fast,\fast,\lambda\}} + (\lambda-1)\mathbf{I} ) \mathbf{1} \right) = \frac{M^2}{6}$ \hfill\refstepcounter{equation}\thetag{\theequation} \label{eq:order_slow_3g}\\
     \botrule
    \end{tabular}
\end{minipage}
\end{subequations}

\begin{subequations}\label{eq:order_fast}
\noindent\begin{minipage}{\textwidth}
    \captionof{table}{Fast order conditions for the MGARK scheme \eqref{eq:MGARK}-\eqref{eq:MGARK_Butcher}}
    \label{tab:MGARK_order_fast}
    \begin{tabular}{p{.05\textwidth} p{.9\textwidth}}
        \toprule
        $p$ & Order condition \\
        \midrule
        1 & $\sum\nolimits_{\lambda=1}^M \mathbf{\bar{b}}^{\{\fast,\lambda\}\top} \mathbf{1} = M$ \hfill\refstepcounter{equation}\thetag{\theequation}\hspace*{.5ex}\label{eq:order_fast_1a} \\
        2 & $\sum\nolimits_{\lambda=1}^M \mathbf{\bar{b}}^{\{\fast,\lambda\}\top} \mathbf{\bar{A}}^{\{\fast,\fast,\lambda\}}\mathbf{1} = \tfrac{M}{2}$ \hfill\refstepcounter{equation}\thetag{\theequation}\hspace*{.5ex}\label{eq:order_fast_2a}\\
        & $\left(\sum\nolimits_{\lambda=1}^M \mathbf{\bar{b}}^{\{\fast,\lambda\}\top} \mathbf{\bar{A}}^{\{\fast,\slow,\lambda\}}\right)\mathbf{1} = \frac{M}{2}$ \hfill\refstepcounter{equation}\thetag{\theequation} \label{eq:order_fast_2b} \\
        3 & $\sum\nolimits_{\lambda=1}^M \mathbf{\bar{b}}^{\{\fast,\lambda\}\top} \diag\left(\mathbf{\bar{A}}^{\{\fast,\fast,\lambda\}}\mathbf{1}\right)\mathbf{\bar{A}}^{\{\fast,\fast,\lambda\}}\mathbf{1} = \frac{M}{3}$ \hfill\refstepcounter{equation}\thetag{\theequation}\hspace*{.5ex}\label{eq:order_fast_3a} \\
        & $\left(\sum\nolimits_{\lambda=1}^M \mathbf{\bar{b}}^{\{\fast,\lambda\}\top}
        \left\{\diag\left(\mathbf{\bar{A}}^{\{\fast,\fast,\lambda\}}\mathbf{1}\right) + (\lambda - 1) \mathbf{I}\right\} \mathbf{\bar{A}}^{\{\fast,\slow,\lambda\}}\right) \mathbf{1} =  \frac{M^2}{3}$ \hfill\refstepcounter{equation}\thetag{\theequation} \label{eq:order_fast_3b} \\
         & $\left( \sum_{\lambda=1}^M \mathbf{\bar{b}}^{\{\fast,\lambda\}\top} \diag\left( \mathbf{\bar{A}}^{\{\fast,\slow,\lambda\}}\mathbf{1}\right)  \mathbf{\bar{A}}^{\{\fast,\slow,\lambda\}}\right)\mathbf{1} = \frac{M}{3}$ \hfill\refstepcounter{equation}\thetag{\theequation} \label{eq:order_fast_3c}\\
         & $\sum\nolimits_{\lambda=1}^M \mathbf{\bar{b}}^{\{\fast,\lambda\}\top} \mathbf{\bar{A}}^{\{\fast,\fast,\lambda\}} \mathbf{\bar{A}}^{\{\fast,\fast,\lambda\}} \mathbf{1} = \frac{M}{6}$ \hfill\refstepcounter{equation}\thetag{\theequation} \label{eq:order_fast_3d}\\
         & $\left( \sum\nolimits_{\lambda=1}^M \mathbf{\bar{b}}^{\{\fast,\lambda\}\top} \!\left\{\mathbf{\bar{A}}^{\{\fast,\fast,\lambda\}} + \left( \sum\nolimits_{\mu=\lambda+1}^M \mathbf{\bar{b}}^{\{\fast,\mu\}\top} \mathbf{1}\!\right)  \mathbf{I}\!\right\}\! \mathbf{\bar{A}}^{\{\fast,\slow,\lambda\}} \right) \!\mathbf{1} \!=\! \frac{M^2}{6},$ \hfill\refstepcounter{equation}\thetag{\theequation} \label{eq:order_fast_3e} \\
         & $\sum_{\lambda=1}^M \mathbf{\bar{b}}^{\{\fast,\lambda\}\top}  \mathbf{\bar{A}}^{\{\fast,\slow,\lambda\}} \left(\sum_{\mu=1}^M \mathbf{\bar{A}}^{\{\slow,\fast,\mu\}}\mathbf{1}\right) = \frac{M^2}{6}$ \hfill\refstepcounter{equation}\thetag{\theequation} \label{eq:order_fast_3f} \\
         & $ \left(\sum_{\lambda=1}^M \mathbf{\bar{b}}^{\{\fast,\lambda\}\top} \mathbf{\bar{A}}^{\{\fast,\slow,\lambda\}}\right)\mathbf{\bar{A}}^{\{\slow,\slow\}}\mathbf{1} = \frac{M}{6}$ \hfill\refstepcounter{equation}\thetag{\theequation} \label{eq:order_fast_3g} \\
         \botrule
    \end{tabular}
\end{minipage}
\end{subequations}

\theorem{\label{theorem:order}For symplectic MGARK schemes, the order-2 conditions {\normalfont(\hyperref[eq:order_slow_2a]{22b})} and {\normalfont(\hyperref[eq:order_fast_2a]{23b})} are redundant. Moreover the following order conditions are equivalent: 
\begin{itemize}
    \item order 2: {\normalfont(\hyperref[eq:order_slow_2b]{22c})} $\Leftrightarrow$ {\normalfont(\hyperref[eq:order_fast_2b]{23c})};
    \item order 3: {\normalfont(\hyperref[eq:order_slow_3a]{22d})} $\Leftrightarrow$ {\normalfont(\hyperref[eq:order_slow_3d]{22g})}, {\normalfont(\hyperref[eq:order_fast_3a]{23d})} $\Leftrightarrow$ {\normalfont(\hyperref[eq:order_fast_3d]{23g})}, {\normalfont(\hyperref[eq:order_slow_3b]{22e})} $\Leftrightarrow$ {\normalfont(\hyperref[eq:order_slow_3e]{22h})} $\Leftrightarrow$ {\normalfont(\hyperref[eq:order_fast_3g]{23j})}, \\ {\normalfont(\hyperref[eq:order_slow_3f]{22i})} $\Leftrightarrow$ {\normalfont(\hyperref[eq:order_fast_3c]{23f})}, {\normalfont(\hyperref[eq:order_slow_3c]{22f})} $\Leftrightarrow$ {\normalfont(\hyperref[eq:order_fast_3f]{23i})}, {\normalfont(\hyperref[eq:order_slow_3g]{22j})} $\Leftrightarrow$ {\normalfont(\hyperref[eq:order_fast_3b]{23e})} $\Leftrightarrow$ {\normalfont(\hyperref[eq:order_fast_3e]{23h})}.
\end{itemize} 
Hence the number of order conditions is reduced due to redundancy as follows:
\begin{itemize}
    \item order 2: the number of 4 order conditions reduces to only 1 order condition;
    \item order 3: the number of 14 order conditions reduces to only 6 order conditions.
\end{itemize}
}\proof{
The MGARK scheme is symplectic and thus satisfies the symplecticity conditions \eqref{eq:symplectic_MGARK_a}-\eqref{eq:symplectic_MGARK_b} of the base schemes. Left multiplication with $\mathbf{1}^\top$ and right multiplication with $\mathbf{1}$ results in the scalar equation
\begin{align*}
    \mathbf{\bar{b}}^{\{\slow\}\top} \mathbf{\bar{A}}^{\{\slow,\slow\}} \mathbf{1} = \tfrac{1}{2}
\intertext{for the slow part. Similarly, summing up all $M$ micro-steps yields}
\sum\limits_{\lambda=1}^M \mathbf{\bar{b}}^{\{\fast,\lambda\}\top} \mathbf{\bar{A}}^{\{\fast,\fast,\lambda\}} \mathbf{1} = \tfrac{M}{2}
\end{align*}
for the fast part.
These are exactly the order-2 conditions {\normalfont(\hyperref[eq:order_slow_2a]{22b})} and {\normalfont(\hyperref[eq:order_fast_2a]{23b})}. W.l.o.g. assume that {\normalfont(\hyperref[eq:order_slow_2b]{22c})} is satisfied. Analogously, algebraic manipulations of the symplecticity condition \eqref{eq:symplectic_MGARK_c} result in the order-2 condition {\normalfont(\hyperref[eq:order_fast_2b]{23c})}. 
According to \cite{gunther2021symplectic}, it holds 
\begin{align*}
    \mathbf{b}^{\{q\}\top} \left(\mathbf{c}^{\{q,m\}} \times \mathbf{c}^{\{q,\ell\}}\right)=\tfrac{1}{3} \ \Leftrightarrow\ \mathbf{b}^{\{m\}\top} \mathbf{A}^{\{m,q\}} \mathbf{c}^{\{q,\ell\}} = \tfrac{1}{6}, \quad q,m,\ell \in \{\slow,\fast\},
\end{align*}
for symplectic GARK schemes. Additionally, as the element-wise product $\times$ commutes, it holds 
\begin{align*}
    \mathbf{b}^{\{q\}\top} \left(\mathbf{c}^{\{q,m\}} \times \mathbf{c}^{\{q,\ell\}}\right)=\tfrac{1}{3} \ \Leftrightarrow\ \mathbf{b}^{\{q\}\top} \left(\mathbf{c}^{\{q,\ell\}} \times \mathbf{c}^{\{q,m\}}\right) = \tfrac{1}{3}, \quad q,m,\ell \in \{\slow,\fast\}.
\end{align*} 
Applying both to the order-3 conditions of the MGARK scheme \eqref{eq:MGARK} results in the desired equivalences of the order-3 conditions. \hfill \endproof}

For an efficient computational process of the sMGARK scheme, the macro-step and the $M$ micro-steps need to stay mostly decoupled. However, the symplecticity conditions \eqref{eq:symplectic_MGARK} introduce a lower bound on the number of zeros in the coupling matrices $\mathbf{\bar{A}}^{\{\slow,\fast,\lambda\}}$ and $\mathbf{\bar{A}}^{\{\fast,\slow,\lambda\}}$.

\lemma{\label{lemma:coupling}For symplectic MGARK schemes it holds 
\begin{align*}\bar{a}_{i,j}^{\{\slow,\fast,\lambda\}} \neq 0 \quad &\text{or} \quad \bar{a}_{j,i}^{\{\fast,\slow,\lambda\}} \neq 0,
\intertext{as well as the implications}
    \bar{a}_{i,j}^{\{\slow,\fast,\lambda\}} = 0 \quad \Rightarrow \quad \bar{a}_{j,i}^{\{\fast,\slow,\lambda\}} = \bar{b}_i^{\{\slow\}}, \quad &\text{and} \quad
    \bar{a}_{j,i}^{\{\fast,\slow,\lambda\}} = 0 \quad \Rightarrow \quad \bar{a}_{i,j}^{\{\slow,\fast,\lambda\}} = \bar{b}_j^{\{\fast,\lambda\}},
\end{align*}
for all $\lambda=1,\ldots,M,\  i=1,\ldots,\bar{s}^{\{\slow\}},\ j=1,\ldots,\bar{s}^{\{\fast,\lambda\}}$.
}
\proof{Consider the symplecticity condition \eqref{eq:symplectic_MGARK_c} of the coupling parts. Assume $\bar{a}_{i,j}^{\{\slow,\fast,\lambda\}} = 0$. Then, the symplecticity condition reduces to $$\bar{a}_{j,i}^{\{\fast,\slow,\lambda\}} \bar{b}_j^{\{\fast,\lambda\}} = \bar{b}_i^{\{\slow\}} \bar{b}_j^{\{\fast,\lambda\}} \quad \Leftrightarrow \quad \bar{a}_{j,i}^{\{\fast,\slow,\lambda\}} = \bar{b}_i^{\{\slow\}}.$$
For $\bar{a}_{j,i}^{\{\fast,\slow,\lambda\}} = 0$, the statement follows analogously. \hfill \endproof}\normalfont 

\remark[decoupled sMGARK schemes]{For decoupled MGARK schemes \cite{sarshar2019design}, it especially has to hold 
\begin{align}\label{eq:decoupled_MGARK}
    \mathbf{\bar{A}}^{\{\slow,\fast,\lambda\}} \times \mathbf{\bar{A}}^{\{\fast,\slow,\lambda\}\top} = \mathbf{0},
\end{align}
where $\times$ denotes element-wise multiplication. In this special case Lemma \ref{lemma:coupling} implies that the coupling matrices can only have entries 
\[
\bar{a}_{i,j}^{\{\slow,\fast,\lambda\}} \in \left\{0,\bar{b}_j^{\{\fast,\lambda\}}\right\}, \quad \bar{a}_{i,j}^{\{\fast,\slow,\lambda\}} \in \left\{0,\bar{b}_j^{\{\slow\}}\right\}.
\]
}\normalfont 

\example[Multirate IMIM2 scheme]{Based on the IMIM2 scheme \cite[Example 3]{gunther2021symplectic}, which is a symplectic and time-reversible GARK scheme, one obtains the multirate IMIM2 (MR-IMIM2) scheme by defining
\begin{align}\label{eq:MR-IMIM2}
    &\mathbf{\bar{A}}^{\{\slow,\slow\}} = \begin{bmatrix} \frac{1}{4} & \beta \\ \frac{1}{2} - \beta & \frac{1}{4}
    \end{bmatrix}, & &\mathbf{\bar{b}}^{\{\slow\}\top} = \begin{bmatrix}
        \tfrac{1}{2} & \tfrac{1}{2}
    \end{bmatrix},  & &\beta \in \mathbb{R}, \nonumber \\
    &\mathbf{\bar{A}}^{\{\fast,\fast,\lambda\}} = \begin{bmatrix} \frac{1}{4} & \alpha \\ \frac{1}{2} - \alpha & \frac{1}{4}
    \end{bmatrix}, &  &\mathbf{\bar{b}}^{\{\fast,\lambda\}\top} = \begin{bmatrix}
        \tfrac{1}{2} & \tfrac{1}{2}
    \end{bmatrix},  & &\alpha \in \mathbb{R},\ \forall \lambda, \\
    &\mathbf{\bar{A}}^{\{\slow,\fast,\lambda\}} = \begin{bmatrix} 0 & 0 \\[0.2ex] \frac{1}{2} & \frac{1}{2}
    \end{bmatrix}, & &\mathbf{\bar{A}}^{\{\fast,\slow,\lambda\}} = \begin{bmatrix} \tfrac{1}{2} & 0 \\[0.2ex] \tfrac{1}{2} & 0
    \end{bmatrix}, & &\forall \lambda. \nonumber
\end{align}
The MR-IMIM2 scheme defines an implicit-implicit (IMIM) symplectic and symmetric MGARK method of order $p=2$ which is algebraically stable.
}\\\normalfont

For $\beta = 0$, the MR-IMIM2 scheme allows for an efficient computational process as one can start with computing $\mathbf{Y}_1^{\{\slow\}}$ (nonlinear system of dimension $n$), then computes the fast stage vectors $\mathbf{Y}_i^{\{\fast,\lambda\}}$ ($M$ nonlinear systems of dimension $2n$), and finally computes $\mathbf{Y}_2^{\{\slow\}}$ (nonlinear system of dimension $n$). For $\alpha = 0$, the computation of the fast stage vectors can be split into $2M$ nonlinear systems of dimension $n$. As $M$ micro-steps of step size $h$ of the IMIM2 scheme would demand the solution of $4M$ ($\alpha = \beta = 0$) or $3M$ ($\beta=0,\ \alpha \neq 0$) nonlinear systems, the MR-IMIM2 scheme demands the solution of less nonlinear systems. 

Furthermore, the MR-IMIM2 scheme evaluates the slow (and often expensive) part $\mathbf{f}^{\{\slow\}}$ for the solution of only $2$ nonlinear systems, whereas the singlerate IMIM2 scheme needs evaluations of $\mathbf{f}^{\{\slow\}}$ for the solution of $2M$ nonlinear systems. Hence the multirate version is expected to be much more efficient for $M \gg 1$.
If $\mathbf{f}^{\{\slow\}}$ is expensive to evaluate, a sMGARK scheme that only needs to evaluate the slow part for the solution of a single nonlinear system may lead to an even more efficient computational process. This is achieved in the following example.

\example[Fastest-first midpoint]{Let $M \in \mathbb{N}$ be even. Consider a sMGARK scheme with base schemes given by the implicit midpoint rule, i.e., 
\[
\left(\mathbf{\bar{A}}^{\{\slow,\slow\}},\mathbf{\bar{b}}^{\{\slow\}} \right) = \left(\mathbf{\bar{A}}^{\{\fast,\fast,\lambda\}},\mathbf{\bar{b}}^{\{\fast,\lambda\}} \right) = \left( \begin{bmatrix}
    \tfrac{1}{2}
\end{bmatrix}, \ \begin{bmatrix}
    1
\end{bmatrix}\right),\quad \lambda=1,\ldots,M.
\]
and the coupling defined via $$\mathbf{\bar{A}}^{\{\slow,\fast,\lambda\}}= \begin{cases} 1, & \lambda = 1,\ldots,M/2, \\ 0, & \lambda = M/2 + 1,\ldots,M, \end{cases} \quad \mathbf{\bar{A}}^{\{\fast,\slow,\lambda\}}= \begin{cases} 0, & \lambda = 1,\ldots,M/2, \\ 1, & \lambda = M/2 + 1,\ldots,M. \end{cases}$$ This scheme is a symplectic and time-reversible MGARK method of order $p=2$ which is algebraically stable.}\normalfont

\remark[Backward error analysis]{Performing a backward error analysis, it can be shown that symplectic schemes preserve a nearby shadow Hamiltonian $\tilde{\H}$ exactly. Based on symplectic NB-series \cite{araujo1997symplectic}, expressions of the shadow Hamiltonian for symplectic GARK schemes have been derived \cite{gunther2021symplectic}. Since MGARK schemes define a subclass of GARK schemes, the backward error analysis can be carried over to sMGARK schemes in a straightforward manner.}\normalfont

\section{Partitioned sMGARK methods for separable Hamiltonians}
\label{sec:P-sMGARK}
In many applications, separable Hamiltonian systems
    $\H(\p,\q) = \mathcal{T}(\p) + \mathcal{V}(\q)$
with kinetic part $\mathcal{T}(\p)$ and potential part $\mathcal{V}(\q)$ occur. Assume that both the kinetic and potential part are split into a slow/expensive and fast/cheap part. Hence, we consider the separable Hamiltonian
\begin{subequations}\label{eq:separable_splitting}
\begin{align}
    \H(\p,\q) &= \mathcal{T}^{\{\slow\}}(\p) + \mathcal{T}^{\{\fast\}}(\p) + \mathcal{V}^{\{\slow\}}(\q) + \mathcal{V}^{\{\fast\}}(\q), \label{eq:MR_separable_Hamiltonian}
\intertext{with four-way partitioned equations of motion}
    \begin{pmatrix}
        \dot{\p} \\ \dot{\q} 
    \end{pmatrix} &= \begin{pmatrix}
        \mathbf{0} \\ \mathcal{T}_\p^{\{\slow\}}(\p)
    \end{pmatrix} + \begin{pmatrix}
        \mathbf{0} \\ \mathcal{T}_\p^{\{\fast\}}(\p)
    \end{pmatrix} + \begin{pmatrix}
        -\mathcal{V}_\q^{\{\slow\}}(\q) \\ \mathbf{0}
    \end{pmatrix} + \begin{pmatrix}
        -\mathcal{V}_\q^{\{\fast\}}(\q) \\ \mathbf{0}
    \end{pmatrix}, \label{eq:sep_H_EoM}
\end{align}
\end{subequations}
where we use abbreviations $\mathcal{T}_\p^{\{q\}} \coloneqq (\nabla_\p \mathcal{T}^{\{q\}})^\top,\ \mathcal{V}_\q^{\{q\}} \coloneqq (\nabla_\q \mathcal{V}^{\{q\}})^\top,\ q \in \{\slow,\fast\}$.
The idea is to solve both slow partitions $\mathcal{T}^{\{\slow\}}(\p)$, $\mathcal{V}^{\{\slow\}}(\q)$ using a large macro-step size $H$, and both fast partitions $\mathcal{T}^{\{\fast\}}(\p)$, $\mathcal{V}^{\{\fast\}}(\q)$ using a small micro-step size $h=H/M$ with a constant multirate factor $M \in \mathbb{N},\ M \gg 1$. One step of a partitioned MGARK scheme, applied to the system \eqref{eq:sep_H_EoM}, advances the solution $(\p_0,\q_0)$ at $t_0$ to the solution $(\p_1,\q_1)$ at time point $t_1 = t_0 + H$ with $M$ equal micro-steps of size $h=H/M$ as follows:
\begin{subequations}\label{eq:pMGARK}
\allowdisplaybreaks
\begin{align}
    \p_1 = \p_0 &- h \sum\limits_{\lambda=1}^{M}\sum\limits_{i=1}^{\bar{s}^{\{\fast,\lambda\}}} \tilde{b}_i^{\{\fast,\lambda\}} \mathcal{V}_\q^{\{\fast\}}\!\left(\mathbf{Q}_i^{\{\fast,\lambda\}}\right) - H \sum\limits_{i=1}^{\bar{s}^{\{\slow\}}} \tilde{b}_i^{\{\slow\}} \mathcal{V}_\q^{\{\slow\}}\!\left(\mathbf{Q}_i^{\{\slow\}}\right), \\
    \q_1 = \q_0 &+ h \sum\limits_{\lambda=1}^{M}\sum\limits_{i=1}^{\tilde{s}^{\{\fast,\lambda\}}} \bar{b}_i^{\{\fast,\lambda\}} \mathcal{T}_\p^{\{\fast\}}\left(\mathbf{P}_i^{\{\fast,\lambda\}}\right) + H \sum\limits_{i=1}^{\tilde{s}^{\{\slow\}}} \bar{b}_i^{\{\slow\}} \mathcal{T}_\p^{\{\slow\}}\left(\mathbf{P}_i^{\{\slow\}}\right), \\
    \begin{split}
    \mathbf{P}_i^{\{\slow\}} = \p_0 &- h\sum\limits_{\lambda=1}^M \sum\limits_{j=1}^{\bar{s}^{\{\fast,\lambda\}}} \tilde{a}_{i,j}^{\{\slow,\fast,\lambda\}} \mathcal{V}_\q^{\{\fast\}}\!\left(\mathbf{Q}_j^{\{\fast,\lambda\}}\right) \\
    &- H \sum\limits_{j=1}^{\bar{s}^{\{\slow\}}} \tilde{a}_{i,j}^{\{\slow,\slow\}} \mathcal{V}_\q^{\{\slow\}}\!\left(\mathbf{Q}_j^{\{\slow\}}\right), \  i=1,\ldots,\tilde{s}^{\{\slow\}}, 
    \end{split}\\
    \begin{split}
    \mathbf{P}_i^{\{\fast,\lambda\}} = \p_0 &- h\sum\limits_{\ell=1}^{\lambda-1} \sum\limits_{j=1}^{\bar{s}^{\{\fast,\ell\}}} \tilde{b}_{j}^{\{\fast,\ell\}} \mathcal{V}_\q^{\{\fast\}}\left(\mathbf{Q}_j^{\{\fast,\ell\}}\right) - H \sum\limits_{j=1}^{\bar{s}^{\{\slow\}}} \tilde{a}_{i,j}^{\{\fast,\slow,\lambda\}} \mathcal{V}_\q^{\{\slow\}}\left(\mathbf{Q}_j^{\{\slow\}}\right) \\
    &- h \sum\limits_{j=1}^{\bar{s}^{\{\fast,\lambda\}}} \tilde{a}_{i,j}^{\{\fast,\fast,\lambda\}} \mathcal{V}_\q^{\{\fast\}}\left(\mathbf{Q}_j^{\{\fast,\lambda\}}\right), \  i = 1,\ldots, \tilde{s}^{\{\fast,\lambda\}}, \, \lambda = 1,\ldots,M, 
    \end{split} \\
    \begin{split}
    \mathbf{Q}_i^{\{\slow\}} = \q_0 &+ h\sum\limits_{\lambda=1}^M \sum\limits_{j=1}^{\tilde{s}^{\{\fast,\lambda\}}} \bar{a}_{i,j}^{\{\slow,\fast,\lambda\}} \mathcal{T}_\p^{\{\fast\}}\!\left(\mathbf{P}_j^{\{\fast,\lambda\}}\right) \\
    &+ H \sum\limits_{j=1}^{\tilde{s}^{\{\slow\}}} \bar{a}_{i,j}^{\{\slow,\slow\}} \mathcal{T}_\p^{\{\slow\}}\!\left(\mathbf{P}_j^{\{\slow\}}\right), \  i=1,\ldots,\bar{s}^{\{\slow\}}, 
    \end{split}\\
    \begin{split}
    \mathbf{Q}_i^{\{\fast,\lambda\}} = \q_0 &+ h\sum\limits_{\ell=1}^{\lambda-1} \sum\limits_{j=1}^{\tilde{s}^{\{\fast,\ell\}}} \bar{b}_{j}^{\{\fast,\ell\}} \mathcal{T}_\p^{\{\fast\}}\left(\mathbf{P}_j^{\{\fast,\ell\}}\right) + H \sum\limits_{j=1}^{\tilde{s}^{\{\slow\}}} \bar{a}_{i,j}^{\{\fast,\slow,\lambda\}} \mathcal{T}_\p^{\{\slow\}}\left(\mathbf{P}_j^{\{\slow\}}\right) \\
    &+ h \sum\limits_{j=1}^{\tilde{s}^{\{\fast,\lambda\}}} \bar{a}_{i,j}^{\{\fast,\fast,\lambda\}} \mathcal{T}_\p^{\{\fast\}}\left(\mathbf{P}_j^{\{\fast,\lambda\}}\right), \  i = 1,\ldots, \bar{s}^{\{\fast,\lambda\}}, \, \lambda = 1,\ldots,M, 
    \end{split}
\end{align}
\end{subequations}
The partitioned MGARK scheme \eqref{eq:pMGARK} has the generalized Butcher tableau 
\begin{align}\label{eq:MGARK_separableHamiltonian_Butcher}
    \raisebox{-18.75pt}{\begin{tabular}{c|c}
         $\mathbf{0}$ & $\mathbf{\hat{A}}$ \\\hline
         $\mathbf{A}$ & $\mathbf{0}$ \\\Xhline{2\arrayrulewidth}\rule{0pt}{3ex} 
         $\mathbf{b}^\top$ & $\mathbf{\hat{b}}^\top$
    \end{tabular}} \raisebox{-23.25pt}{$~\coloneqq~$}
    \begin{tabular}{c c|c c}
        $\ \mathbf{0}$& $\mathbf{0}$ & $\mathbf{\hat{A}}^{\{\fast,\fast\}}$ & $\mathbf{\hat{A}}^{\{\fast,\slow\}}$ \\\rule{0pt}{3ex}
        $\mathbf{0}$ & $\mathbf{0}$ & $\mathbf{\hat{A}}^{\{\slow,\fast\}}$ & $\mathbf{\hat{A}}^{\{\slow,\slow\}}$ \\\hline\rule{0pt}{3ex}
         $\mathbf{A}^{\{\fast,\fast\}}$ & $\mathbf{A}^{\{\fast,\slow\}}$ & $\mathbf{0}$ & $\mathbf{0}$ \\\rule{0pt}{3ex}
          $\mathbf{A}^{\{\slow,\fast\}}$ & $\mathbf{A}^{\{\slow,\slow\}}$ & $\mathbf{0}$ & $\mathbf{0}$ \\ \Xhline{2\arrayrulewidth}\rule{0pt}{3ex}
          $\mathbf{b}^{\{\fast\}\top}$ & $\mathbf{b}^{\{\slow\}\top}$ & $\mathbf{\hat{b}}^{\{\fast\}\top}$ & $\mathbf{\hat{b}}^{\{\slow\}\top}$
    \end{tabular} 
\end{align}
with left half (ignoring the zero-blocks) as in \eqref{eq:MGARK_Butcher} and right half (ignoring the zero-blocks)
\begin{align}\label{eq:pMGARK_Butcher-Hat}
    \raisebox{-28.55pt}{\begin{tabular}{c|c}
         $\mathbf{\hat{A}}^{\{\fast,\fast\}}$ & $\mathbf{\hat{A}}^{\{\fast,\slow\}}$ \\\hline\rule{0pt}{3ex}
         $\mathbf{\hat{A}}^{\{\slow,\fast\}}$ & $\mathbf{\hat{A}}^{\{\slow,\slow\}}$ \\ \Xhline{2\arrayrulewidth}\rule{0pt}{3ex}
         $\mathbf{\hat{b}}^{\{\fast\}\top}$ & $\mathbf{\hat{b}}^{\{\slow\}\top}$
    \end{tabular}} \!\raisebox{-35.375pt}{$~\coloneqq~$}\! \begin{tabular}{c c c c | c}
         $\frac{1}{M} \mathbf{\tilde{A}}^{\{\fast,\fast,1\}}$ & $\mathbf{0}$ & $\cdots$ & $\mathbf{0}$ & $\mathbf{\tilde{A}}^{\{\fast,\slow,1\}}$ \\
         $\frac{1}{M} \mathbf{1} \mathbf{\tilde{b}}^{\{\fast,1\}\top}$ & $\frac{1}{M} \mathbf{\tilde{A}}^{\{\fast,\fast,2\}}$ & $\ddots$ & $\vdots$ & $\mathbf{\tilde{A}}^{\{\fast,\slow,2\}}$ \\
         $\vdots$ & $\ddots$ & $\ddots$ & $\mathbf{0}$ & $\vdots$ \\
         $\frac{1}{M} \mathbf{1} \mathbf{\tilde{b}}^{\{\fast,1\}\top}$ & $\cdots$ & $\frac{1}{M} \mathbf{1} \mathbf{\tilde{b}}^{\{\fast,M-1\}\top}$ & $\frac{1}{M} \mathbf{\tilde{A}}^{\{\fast,\fast,M\}}$ & $\mathbf{\tilde{A}}^{\{\fast,\slow,M\}}$ \\[-2ex]
         & & & & 
         \\\hline\rule{0pt}{3ex}
         $\frac{1}{M} \mathbf{\tilde{A}}^{\{\slow,\fast,1\}}$ & $\frac{1}{M} \mathbf{\tilde{A}}^{\{\slow,\fast,2\}}$ & $\cdots$ & $\frac{1}{M} \mathbf{\tilde{A}}^{\{\slow,\fast,M\}}$ & $\mathbf{\tilde{A}}^{\{\slow,\slow\}}$ \\[-2ex]
         & & & & \\  \Xhline{2\arrayrulewidth}\rule{0pt}{3ex}
         $\frac{1}{M} \mathbf{\tilde{b}}^{\{\fast,1\}\top}$ & $\frac{1}{M} \mathbf{\tilde{b}}^{\{\fast,2\}\top}$ & $\cdots$ & $\frac{1}{M} \mathbf{\tilde{b}}^{\{\fast,M\}\top}$ & $\mathbf{\tilde{b}}^{\{\slow\}\top}$
    \end{tabular}\,.
\end{align}
In this case one is interested in sMGARK methods for four-way partitioned Hamiltonian systems where the separability of the Hamiltonian system introduces a special structure given by the zero-blocks in \eqref{eq:MGARK_separableHamiltonian_Butcher}. This special structure results in simplified conditions on symplecticity and symmetry.

\begin{subequations}
\definition[Symmetric partitioned MGARK scheme]{The partitioned MGARK scheme \eqref{eq:pMGARK}-\eqref{eq:MGARK_separableHamiltonian_Butcher} is symmetric if both blocks \eqref{eq:MGARK_Butcher} and \eqref{eq:pMGARK_Butcher-Hat} satisfy the symmetry conditions \eqref{eq:symmetric_MGARK} for MGARK schemes.
}\normalfont
\end{subequations}

\definition[Symplectic partitioned MGARK scheme]{The partitioned MGARK scheme \eqref{eq:pMGARK}-\eqref{eq:MGARK_separableHamiltonian_Butcher} is symplectic if
\begin{subequations}
    \begin{align}
        \mathbf{\bar{A}}^{\{\slow,\slow\}\top} \mathbf{\tilde{B}}^{\{\slow\}} + \mathbf{\bar{B}}^{\{\slow\}} \mathbf{\tilde{A}}^{\{\slow,\slow\}} - \mathbf{\bar{b}}^{\{\slow\}} \mathbf{\tilde{b}}^{\{\slow\}\top} &= \mathbf{0}, \label{eq:symplectic_pMGARK_a}
    \intertext{and for all $\lambda=1,\ldots,M$ it holds}
        \mathbf{\bar{A}}^{\{\fast,\fast,\lambda\}\top} \mathbf{\tilde{B}}^{\{\fast,\lambda\}} + \mathbf{\bar{B}}^{\{\fast,\lambda\}} \mathbf{\tilde{A}}^{\{\fast,\fast,\lambda\}} - \mathbf{\bar{b}}^{\{\fast,\lambda\}} \mathbf{\tilde{b}}^{\{\fast,\lambda\}\top} &= \mathbf{0}, \label{eq:symplectic_pMGARK_b} \\ 
        \mathbf{\bar{A}}^{\{\fast,\slow,\lambda\}\top} \mathbf{\tilde{B}}^{\{\fast,\lambda\}} + \mathbf{\bar{B}}^{\{\slow\}} \mathbf{\tilde{A}}^{\{\slow,\fast,\lambda\}} - \mathbf{\bar{b}}^{\{\slow\}} \mathbf{\tilde{b}}^{\{\fast,\lambda\}\top} &= \mathbf{0}, \label{eq:symplectic_pMGARK_c} \\
        \mathbf{\bar{A}}^{\{\slow,\fast,\lambda\}\top} \mathbf{\tilde{B}}^{\{\slow\}} + \mathbf{\bar{B}}^{\{\fast,\lambda\}} \mathbf{\tilde{A}}^{\{\fast,\slow,\lambda\}} - \mathbf{\bar{b}}^{\{\fast,\lambda\}} \mathbf{\tilde{b}}^{\{\slow\}\top} &= \mathbf{0}. \label{eq:symplectic_pMGARK_d}
    \end{align}
\end{subequations}
}\normalfont

The order conditions up to third order look similar to those in Tables \ref{tab:MGARK_order_slow} and \ref{tab:MGARK_order_fast}. The number of conditions doubles analogously to partitioned Runge--Kutta schemes. For completeness, the order conditions for partitioned MGARK schemes can be found in Appendix \ref{sec:order-PMGARK}. Analogously to the results in Section \ref{sec:sMGARK}, the order-two conditions vanish for symmetric schemes. For symplectic partitioned MGARK schemes, the order conditions reduce similar to the results in Theorem \ref{theorem:order}. Specifically, the number of $8$ order-two conditions reduces to only $2$ order conditions, and the number of $28$ order-three conditions reduces to only $12$ order conditions.

\remark{Setting $\mathbf{A}^{\{q,m\}} = \mathbf{\hat{A}}^{\{q,m\}}$ and $\mathbf{b}^{\{q\}} = \mathbf{\hat{b}}^{\{q\}}$ for $q,m\in\{\slow,\fast\}$, the partitioned MGARK scheme \eqref{eq:pMGARK} reduces to the MGARK scheme \eqref{eq:MGARK} with $\f^{\{\slow\}\top} = \begin{pmatrix}
    -\mathcal{V}_\q^{\{\slow\}\top} & \mathcal{T}_\p^{\{\slow\}\top} 
\end{pmatrix}$ and $\f^{\{\fast\}\top} = \begin{pmatrix}
    -\mathcal{V}_\q^{\{\fast\}\top} & \mathcal{T}_\p^{\{\fast\}\top} 
\end{pmatrix}$.}\\\normalfont

The particular splitting \eqref{eq:sep_H_EoM} allows for the derivation of explicit partitioned sMGARK schemes. A partitioned sMGARK scheme is explicit if it holds \begin{align}\label{eq:pMGARK_explicit}
    \mathbf{A}^{\{q,m\}} \times \mathbf{\hat{A}}^{\{m,q\}\top} = \mathbf{0}, \ q,m\in \{\slow,\fast\},
\end{align}
with $\times$ denoting element-wise multiplication. In this case, the scheme \eqref{eq:pMGARK} turns into a composition of flows to the subsystems in \eqref{eq:sep_H_EoM}, i.e., it is a multirate scheme within the framework of splitting and composition methods where the flows of the subsystems can be computed exactly. They are given by 
\begin{align}\label{eq:exact_flows}
\varphi_t^{\{q\}}(\y_0) = \y_0 + t \mathbf{f}^{\{q\}}(\y_0),\ q =1,\ldots,4.
\end{align}
One example of an explicit partitioned sMGARK scheme is obtained based on the leapfrog scheme.
 
\begin{subequations}\label{eq:nested_leapfrog}
\example[MR-LPFR]{Let $M \in \mathbb{N}$ be even.
The multirate leapfrog (MR-LPFR) scheme is an explicit, symplectic and symmetric partitioned MGARK scheme \eqref{eq:pMGARK} with generalized Butcher tableau \eqref{eq:MGARK_separableHamiltonian_Butcher} defined via 
\begin{align}
    &\mathbf{\bar{A}}^{\{\slow,\slow\}} = \mathbf{\bar{A}}^{\{\fast,\fast,\lambda\}} = \begin{bmatrix} 0 & 0 \\ \tfrac{1}{2} & \tfrac{1}{2} \end{bmatrix},\ \forall \lambda, &
    &\mathbf{\tilde{A}}^{\{\slow,\slow\}} = \mathbf{\tilde{A}}^{\{\fast,\fast,\lambda\}} = \begin{bmatrix} \tfrac{1}{2} & 0 \\[0.1cm] \tfrac{1}{2} & 0 \end{bmatrix},\ \forall \lambda, \\
    &\mathbf{\bar{A}}^{\{\slow,\fast,\lambda\}} \!=\! \begin{cases} \begin{bmatrix} 0 & 0 \\ \tfrac{1}{2} & \tfrac{1}{2} \end{bmatrix}, & \lambda\!=\!1,\ldots,\tfrac{M}{2}, \\[0.4cm] \begin{bmatrix}0 & 0 \\ 0 & 0\end{bmatrix}, & \lambda \!=\! \tfrac{M+2}{2},\ldots,M,\end{cases} & &\mathbf{\bar{A}}^{\{\fast,\slow,\lambda\}} \!=\! \begin{cases} \begin{bmatrix} 0 & 0 \\ \tfrac{1}{2} & \tfrac{1}{2} \end{bmatrix}, & \lambda=1, \\[0.4cm] \begin{bmatrix}\tfrac{1}{2} & \tfrac{1}{2} \\ \tfrac{1}{2} & \tfrac{1}{2}\end{bmatrix}, & \lambda = 2,\ldots,M,\end{cases} \\ 
    &\mathbf{\tilde{A}}^{\{\fast,\slow,\lambda\}} \!=\! \begin{cases} \begin{bmatrix} \tfrac{1}{2} & 0 \\ \tfrac{1}{2} & 0 \end{bmatrix}, & \lambda\!=\!1,\ldots,\tfrac{M}{2}, \\[0.4cm] \begin{bmatrix}\tfrac{1}{2} & \tfrac{1}{2} \\ \tfrac{1}{2} & \tfrac{1}{2}\end{bmatrix}, & \lambda \!=\! \tfrac{M+2}{2},\ldots,M,\end{cases} & &\mathbf{\tilde{A}}^{\{\slow,\fast,\lambda\}} \!=\! \begin{cases} \begin{bmatrix} \tfrac{1}{2} & 0 \\ \tfrac{1}{2} & 0 \end{bmatrix}, & \lambda=1, \\[0.4cm] \begin{bmatrix}0 & 0 \\ 0 & 0\end{bmatrix}, & \lambda = 2,\ldots,M,\end{cases} \\
    &\mathbf{\bar{b}}^{\{\slow\}\top} = \mathbf{\bar{b}}^{\{\fast\}\top} = \begin{bmatrix} \tfrac{1}{2} & \tfrac{1}{2} \end{bmatrix},\ \forall \lambda, & &\mathbf{\tilde{b}}^{\{\slow\}\top} = \mathbf{\tilde{b}}^{\{\fast\}\top} = \begin{bmatrix} \tfrac{1}{2} & \tfrac{1}{2} \end{bmatrix},\ \forall \lambda.
\end{align}
}\normalfont
\end{subequations}
The MR-LPFR scheme \eqref{eq:nested_leapfrog} satisfies \eqref{eq:pMGARK_explicit}, i.e., it is indeed explicit. A corresponding pseudocode of the scheme can be found in Appendix \ref{sec:MR-leapfrog}. 

On the one hand, applying a partitioned sMGARK scheme \eqref{eq:pMGARK} satisfying \eqref{eq:pMGARK_explicit} to a four-way partitioned system \eqref{eq:separable_splitting} yields an explicit integration scheme. On the other hand, \eqref{eq:pMGARK_explicit} implies that the subsystems are integrated independently. The subsystems in \eqref{eq:sep_H_EoM} are not stable from an analytical point of view as it can be seen by their exact flows \eqref{eq:exact_flows}. This may impose strong step size restrictions if one of the subsystems is stiff.
Often, the fast subsystem is stiff and thus imposes a strong step size restriction on the micro-step size $h$. One may overcome this by choosing the multirate factor $M$ sufficiently large but this lowers the computational efficiency. Choosing an implicit scheme for the stiff subsystem may be beneficial in order to remove the step size restriction on $h$, so that one can choose $M$ only based on the accuracy demands. As it is still desirable to have an explicit integrator for the non-stiff partitions, a combination of implicit and explicit schemes could be very efficient.

\section{Implicit-explicit sMGARK schemes}
\label{sec:IMEX-sMGARK}
In this section, we pay attention to the special case of potential splitting $\mathcal{V}(\q) = \mathcal{V}^{\{\slow\}}(\q) + \mathcal{V}^{\{\fast\}}(\q)$. The slow force $\mathcal{V}_{\q}^{\{\slow\}}(\q)$ is assumed to be non-stiff and expensive to evaluate, whereas the fast force $\mathcal{V}_{\q}^{\{\fast\}}(\q)$ is stiff and characterized by cheap evaluation costs. The kinetic part may contain both fast and slow components, and as it frequently is of the form $\mathcal{T}(\p) = \sum_{i=1}^{n_p} \frac{1}{m_i} p_i^2$ and thus cheap to evaluate, we treat it as a fast partition. Consequently, we end up with the two-way partitioned Hamiltonian system
\begin{align}
    \H(\p,\q) = \H^{\{\fast\}}(\p,\q) + \H^{\{\slow\}}(\q),
\end{align}
with Hamiltonian equations of motion 
\begin{align}\label{eq:IMEX_EoM}
    \begin{pmatrix}
        \dot{\p} \\ \dot{\q}
    \end{pmatrix} &= \begin{pmatrix}
        -\mathcal{V}_\q^{\{\fast\}}(\q) \\ \mathcal{T}_\p(\p)
    \end{pmatrix} + \begin{pmatrix}
        -\mathcal{V}_\q^{\{\slow\}}(\q) \\ \mathbf{0}
    \end{pmatrix}= \f^{\{\slow\}}(\q) + \f^{\{\fast\}}(\p,\q).
\end{align}
As the fast subsystem $\f^{\{\fast\}}(\p,\q)$ is assumed to be stiff, it is beneficial to treat it implicitly. The non-stiff subsystem $\f^{\{\slow\}}(\q)$ is preferably solved explicitly using large time steps in order to save unnecessary evaluations of the expensive force.

\theorem{Consider a sMGARK scheme \eqref{eq:MGARK} that is decoupled, i.e., it satisfies the condition \eqref{eq:decoupled_MGARK}. Then the sMGARK scheme, applied to \eqref{eq:IMEX_EoM}, results in an implicit-explicit (IMEX) integration scheme that is explicit for the slow and expensive subsystem $\f^{\{\slow\}}(\q)$, implicit for the fast subsystem $\f^{\{\fast\}}(\p,\q)$. }
\proof{
Exploiting the special structure of the system \eqref{eq:IMEX_EoM}, the formulae for the slow stage vectors reads 
\begin{align*}
    \begin{pmatrix}
        \mathbf{P}_i^{\{\slow\}} \\ \mathbf{Q}_i^{\{\slow\}}  
    \end{pmatrix} = \begin{pmatrix}
        \p_0 \\ \q_0
    \end{pmatrix} &+ H \sum\limits_{j=1}^{\bar{s}^{\{\slow\}}} \bar{a}_{i,j}^{\{\slow,\slow\}} \begin{pmatrix}
        -\mathcal{V}_\q^{\{\slow\}}\left(\mathbf{Q}_j^{\{\slow\}}\right) \\ \mathbf{0}
    \end{pmatrix} \\
    &+ h \sum\limits_{\lambda = 1}^M \sum\limits_{j=1}^{\bar{s}^{\{\fast,\lambda\}}} \bar{a}_{i,j}^{\{\slow,\fast,\lambda\}} \begin{pmatrix}
        -\mathcal{V}_\q^{\{\fast\}}\left(\mathbf{Q}_j^{\{\fast,\lambda\}} \right) \\ \mathcal{T}_\p\left(\mathbf{P}_j^{\{\fast,\lambda\}}\right)
    \end{pmatrix}.
\end{align*} 
As the MGARK scheme is assumed to be decoupled, the last term can be computed explicitly. The special structure in the $H$-term allows for an explicit computation. Hence the overall computation of the slow stage vectors is explicit. On the other hand, the formula for the fast stage vectors reads
\begin{align*}
    \begin{pmatrix}
        \mathbf{P}_i^{\{\fast,\lambda\}} \\ \mathbf{Q}_i^{\{\fast,\lambda\}}
    \end{pmatrix} = &\begin{pmatrix}
        \p_0 \\ \q_0
    \end{pmatrix} + h \sum\limits_{\ell=1}^{\lambda-1} \sum\limits_{j=1}^{\bar{s}^{\{\fast,\ell\}}} \bar{b}_j^{\{\fast,\ell\}} \begin{pmatrix}
        -\mathcal{V}_\q^{\{\fast\}}\left( \mathbf{Q}_j^{\{\fast,\ell\}}  \right) \\ \mathcal{T}_\p \left(\mathbf{P}_j^{\{\fast,\ell\}}\right) 
    \end{pmatrix} \\
    &+ H \sum \limits_{j=1}^{\bar{s}^{\{\slow\}}} \bar{a}_{i,j}^{\{\fast,\slow,\lambda\}} \begin{pmatrix}
        -\mathcal{V}_\q^{\{\slow\}} \left(\mathbf{Q}_j^{\{\slow\}}\right) \\ \mathbf{0}
    \end{pmatrix} + h \sum\limits_{j=1}^{\bar{s}^{\{\fast,\lambda\}}} \bar{a}_{i,j}^{\{\fast,\fast,\lambda\}} \begin{pmatrix}
        -\mathcal{V}_\q^{\{\fast\}} \left(\mathbf{Q}_j^{\{\fast,\lambda\}}\right) \\ \mathcal{T}_\p\left(\mathbf{P}_j^{\{\fast,\lambda\}}\right)
    \end{pmatrix},
\end{align*}
for $i=1,\ldots,\bar{s}^{\{\fast,\lambda\}}, \quad \lambda = 1,\ldots,M.$
In general, the computation of the fast stage vectors demands the solution of a nonlinear system because the base schemes $(\mathbf{\bar{A}}^{\{\fast,\fast,\lambda\}},\mathbf{\bar{b}}^{\{\fast,\lambda\}})$ define symplectic Runge--Kutta methods that are implicit and algebraically stable. Hence, the resulting method is an IMEX scheme that is explicit for the slow and expensive subsystem $\f^{\{\slow\}}(\q)$ but at the same time is implicit for the fast subsystem $\f^{\{\fast\}}(\p,\q)$. \hfill\endproof}

\remark{Due to the zero-block in the slow subsystem, the computation of the slow stage vectors $\mathbf{P}_i^{\{\slow\}}$ is redundant.}\\\normalfont

Although the IMEX sMGARK scheme is symplectic and thus algebraically stable, step size restrictions on the macro-step size $H$ may occur since the slow subsystem with exact flow $\varphi_t^{\{\slow\}}(\y_0) = \y_0 + t \mathbf{f}^{\{\slow\}}(\y_0)$ is not stable from an analytical point of view. As it is assumed that the slow subsystem is non-stiff, the constraint on $H$ is typically dominated by the accuracy demands.  
All in all, the IMEX approach comes with improved stability properties compared to the fully explicit schemes discussed in Section 3. Furthermore, it saves expensive evaluations of the slow subsystem by using an explicit integration scheme for it, resulting in a more efficient integration process compared to the IMIM schemes discussed in Section \ref{sec:sMGARK}. This makes the IMEX approach preferable in case of separable Hamiltonian systems with multirate behavior and stiff subsystem $\f^{\{\fast\}}(\p,\q)$. One example for an IMEX sMGARK scheme is obtained by combining the ideas of the implicit midpoint rule and the leapfrog method.

\example[MR-IMEX2]{Consider the sMGARK scheme \eqref{eq:MGARK} given by the generalized Butcher tableau 
\begin{align}
    \mathbf{\bar{A}}^{\{\slow,\slow\}} &= \begin{bmatrix}
        \tfrac{1}{4} & 0 \\ \tfrac{1}{2} & \tfrac{1}{4}
    \end{bmatrix}, & \mathbf{\bar{b}}^{\{\slow\}\top} &= \begin{bmatrix} \tfrac{1}{2} &\tfrac{1}{2} \end{bmatrix},  \nonumber \\
    \mathbf{\bar{A}}^{\{\fast,\fast,\lambda\}} &= \begin{bmatrix}
        \tfrac{1}{2}
    \end{bmatrix}, & \mathbf{\bar{b}}^{\{\fast,\lambda\}\top} &= \begin{bmatrix} 1 \end{bmatrix}, \quad \forall \lambda, \label{eq:IMEX2} \\
    \mathbf{\bar{A}}^{\{\slow,\fast,\lambda\}} &= \begin{bmatrix}
        0 \\ 1
    \end{bmatrix}, & \mathbf{\bar{A}}^{\{\fast,\slow,\lambda\}} &= \begin{bmatrix}
        \tfrac{1}{2} & 0
    \end{bmatrix}, \quad \forall \lambda. \nonumber
\end{align}
This is, in general, an implicit-implicit sMGARK scheme of order two that is symmetric and symplectic. However, for partitions of the form \eqref{eq:IMEX_EoM}, this scheme becomes an IMEX scheme.
}\normalfont

The MR-IMEX2 scheme \eqref{eq:IMEX2} is an impulse method, i.e., the slow part $\f^{\{\slow\}}$ is only used at the beginning and at the end of each macro-step to update the momenta $\p$. In between, there are $M$ micro-steps of the implicit midpoint rule, applied to the fast subsystem $\f^{\{\fast\}}$. Moreover, one can significantly reduce the number of evaluations of $\f^{\{\slow\}}$. The second evaluation of $\f^{\{\slow\}}$ in macro-step $n$ and the first evaluation of $\f^{\{\slow\}}$ at macro-step $n+1$ coincide. Hence one can easily combine these two updates so that computing $\tilde{N}$ macro-steps demands $\tilde{N}+1$ evaluations of the expensive right-hand side $\f^{\{\slow\}}$.

\remark{Both examples in Section \ref{sec:sMGARK} also reduce to an IMEX sMGARK scheme of order two when applied to systems \eqref{eq:IMEX_EoM}.}\normalfont

\section{Construction of higher-order schemes based on advanced composition}
\label{sec:composition}
In the previous sections, we have discussed examples of (partitioned) sMGARK of order two. Composition techniques provide a tool for obtaining (partitioned) sMGARK schemes of arbitrarily high convergence order. In this section, we will consider compositions of sMGARK schemes. Compositions of partitioned sMGARK schemes are obtained by applying the composition technique to both blocks \eqref{eq:MGARK_Butcher} and \eqref{eq:pMGARK_Butcher-Hat} in \eqref{eq:MGARK_separableHamiltonian_Butcher}. We will start with the general composition approach \cite{suzuki1990fractal,yoshida1990construction}.

\subsection{Composition schemes} 
Let $\Phi_H$ denote a sMGARK method of order $k$ that is defined via the generalized Butcher tableau \eqref{eq:MGARK_Butcher}. Consider the composition
\begin{align}\label{eq:composition_scheme}
    \Psi_H = \Phi_{\gamma_r H} \circ \ldots \circ \Phi_{\gamma_1 H}
\end{align}
with $\Phi_H$ as the underlying base scheme and weights $\gamma_i$ satisfying the symmetry conditions $\gamma_i = \gamma_{r+1-i},\ i=1,\ldots,r$. The resulting composition scheme is again a symmetric (symmetric composition of symmetric methods) and symplectic (composition of symplectic methods) MGARK scheme with generalized Butcher tableau \cite{gonzalez2022unified}
\begin{align}\label{eq:composition_Butcher}
    \scalebox{.85}{
     \begin{tabular}{c c c c|c c c c}
         $\gamma_1 \mathbf{A}^{\{\fast,\fast\}}$ & $\mathbf{0}$ & $\cdots$ & $\mathbf{0}$ & $\gamma_1 \mathbf{A}^{\{\fast,\slow\}}$ & $\mathbf{0}$ & $\cdots$ & $\mathbf{0}$ \\
         $\gamma_1 \mathbf{1} \mathbf{b}^{\{\fast\}\top}$ & $\gamma_2 \mathbf{A}^{\{\fast,\fast\}}$ & $\ddots$ & $\vdots$ & $\gamma_1 \mathbf{1} \mathbf{b}^{\{\slow\}\top}$ & $\gamma_2 \mathbf{A}^{\{\fast,\slow\}}$ & $\ddots$ & $\vdots$ \\
         $\vdots$ & $\ddots$ & $\ddots$ & $\mathbf{0}$ & $\vdots$ & $\ddots$ & $\ddots$ & $\mathbf{0}$ \\
         $\gamma_1 \mathbf{1} \mathbf{b}^{\{\fast\}\top}$ & $\cdots$ & $\gamma_{r-1} \mathbf{1} \mathbf{b}^{\{\fast\}\top}$ & $\gamma_r \mathbf{A}^{\{\fast,\fast\}}$ & $\gamma_1 \mathbf{1} \mathbf{b}^{\{\slow\}\top}$ & $\cdots$ & $\gamma_{r-1} \mathbf{1} \mathbf{b}^{\{\slow\}\top}$ & $\gamma_r \mathbf{A}^{\{\fast,\slow\}}$\\\hline\rule{0pt}{3ex}
         $\gamma_1 \mathbf{A}^{\{\slow,\fast\}}$ & $\mathbf{0}$ & $\cdots$ & $\mathbf{0}$ & $\gamma_1 \mathbf{A}^{\{\slow,\slow\}}$ & $\mathbf{0}$ & $\cdots$ & $\mathbf{0}$ \\
         $\gamma_1 \mathbf{1} \mathbf{b}^{\{\fast\}\top}$ & $\gamma_2 \mathbf{A}^{\{\slow,\fast\}}$ & $\ddots$ & $\vdots$ & $\gamma_1 \mathbf{1} \mathbf{b}^{\{\slow\}\top}$ & $\gamma_2 \mathbf{A}^{\{\slow,\slow\}}$ & $\ddots$ & $\vdots$ \\
         $\vdots$ & $\ddots$ & $\ddots$ & $\mathbf{0}$ & $\vdots$ & $\ddots$ & $\ddots$ & $\mathbf{0}$ \\
         $\gamma_1 \mathbf{1} \mathbf{b}^{\{\fast\}\top}$ & $\cdots$ & $\gamma_{r-1} \mathbf{1} \mathbf{b}^{\{\fast\}\top}$ & $\gamma_r \mathbf{A}^{\{\slow,\fast\}}$ & $\gamma_1 \mathbf{1} \mathbf{b}^{\{\slow\}\top}$ & $\cdots$ & $\gamma_{r-1} \mathbf{1} \mathbf{b}^{\{\slow\}\top}$ & $\gamma_r \mathbf{A}^{\{\slow,\slow\}}$\\ \Xhline{2\arrayrulewidth}\rule{0pt}{3ex}
         $\gamma_1 \mathbf{b}^{\{\fast\}\top}$ &  $\gamma_2 \mathbf{b}^{\{\fast\}\top}$ & $\cdots$ &  $\gamma_r \mathbf{b}^{\{\fast\}\top}$ & $\gamma_1 \mathbf{b}^{\{\slow\}\top}$ & $\gamma_2 \mathbf{b}^{\{\slow\}\top}$ & $\cdots$ & $\gamma_r \mathbf{b}^{\{\slow\}\top}$
    \end{tabular} .
    }
\end{align}
The composition scheme \eqref{eq:composition_scheme}-\eqref{eq:composition_Butcher} is a MGARK scheme consisting of $r \cdot M$ micro-steps with micro-step sizes $r \cdot \gamma_i \tfrac{H}{rM} = \gamma_i h$. If the weights $\gamma_i$ satisfy the conditions \cite{HairerLubichWanner}
\begin{subequations}\label{eq:composition_conditions}
\begin{align}
    \gamma_1 + \ldots + \gamma_r &= 1, \\
    \gamma_1^{p+1} + \ldots + \gamma_r^{p+1} &= 0,  \label{eq:composition_condition2}
\end{align}
\end{subequations}
then the composition scheme \eqref{eq:composition_scheme}-\eqref{eq:composition_Butcher} has order $p+2$ (as symmetric schemes are always of even order).

\remark{Note that the second condition \eqref{eq:composition_condition2} demands to have some negative $\gamma_i$. Even if the underlying sMGARK scheme $\Phi_H$ is algebraically stable, the stability domain of the composition scheme will show holes in the left half-plane \cite{li1997raising}.}\\\normalfont

Choosing one of the previously introduced sMGARK schemes of order two as the underlying base scheme $\Phi_H$, one is able to derive sMGARK schemes of arbitrarily high order. There are two common approaches to achieve this.
\example[Triple jump]{
The smallest value of $r$ that allows for a real solution of \eqref{eq:composition_conditions} is $r=3$. Imposing symmetry ($\gamma_1 = \gamma_3$), the unique real solution is given by \cite{yoshida1990construction}
\begin{align}\label{eq:triple-jump}
    \gamma_1 = \gamma_3 = \frac{1}{2 - 2^{1/(p+1)}}, \qquad \gamma_2 = - \frac{2^{1/(p+1)}}{2-2^{1/(p+1)}}.
\end{align}
}\normalfont
Here, one has weights $\lvert \gamma_i \rvert > 1$, i.e., we do not stay within the time window $[t_0,t_0 + H]$. This fact may be permissible in orbit calculations, but they might be harmful in situations when true solutions $y(t)$ pass too near singularities \cite{kahan1997composition}. The idea of Suzuki \cite{suzuki1990fractal} is avoiding this at the price of two additional evaluations of the underlying base scheme ($r=5$).

\example[Suzuki's fractals]{For $r=5$, the best solution of \eqref{eq:composition_conditions} is given by \cite{suzuki1990fractal}
\begin{align}\label{eq:suzuki_fractals}
    \gamma_1 = \gamma_2 = \gamma_4 = \gamma_5 = \frac{1}{4 - 4^{1/(p+1)}}, \quad \gamma_3 = -\frac{4^{1/(p+1)}}{4 - 4^{1/(p+1)}}.
\end{align}
}\normalfont
In this case, it holds $\sum_{i=1}^l \gamma_i \in [0,1]$ for all $l=1,\ldots,r$ so that all step points remain in the time window $[t_0,t_0 + H]$. 

Both ideas can be applied in a recursive manner to obtain symplectic, symmetric and multirate-exploiting GARK schemes of arbitrarily high convergence order. However, the recursive application of these approaches is very inefficient as the number of applications of the underlying basis scheme of order two increases drastically.
A more efficient way of directly deriving composition schemes of a certain order $p \leq 10$ is given by advanced composition schemes \cite{omelyan2002construction}. 

\subsection{Advanced composition schemes}

The idea of advanced composition schemes is to directly write down all conditions that have to be satisfied to derive a composition scheme of order $p$. This approach allows for the construction of composition schemes of higher order with way fewer calls to the underlying base scheme. In \cite{omelyan2002construction}, recursive relations that define the nonlinear systems starting from a base scheme of order $k$ to derive symmetric composition schemes up to order $k+8$ are carried out. As we have derived sMGARK schemes of order two, we are particularly interested in the special case $k=2$. This special case has been investigated in \cite{kahan1997composition}. For specific sets of coefficients, we refer to the appendix therein. In Table \ref{tab:composition_steps}, we have summarized the number of applications of the underlying base scheme $\Phi_H$ that is needed in order to obtain a composition scheme $\Psi_H$ of order $p \in \{4,6,8,10\}$ by using the triple jump (TJ), Suzuki's fractals (SF), advanced composition (AC), advanced composition s.t. one stays in the time window $[t_0,t_0+H]$ (AC*). This underlines the advantage of advanced composition, especially if one wants to solve initial value problems where the true solution $y(t)$ pass too near singularities and thus has to stay inside the time window $[t_0, t_0 + H]$. 

\begin{table}[!htp]
    \caption{Comparison of the number of required composition steps $r$ to obtain a certain order $p$ using the triple jump, Suzuki's fractals, and advanced composition schemes.}
    \label{tab:composition_steps}
    \begin{tabular}{|p{3cm}|p{.5cm}p{1cm}p{1cm}p{1cm}p{1cm}|}
    \hline
    order $p$ of $\Psi_H$ &  & 4 & 6 & 8 & 10 \\\hline
    TJ &  & 3 & 9 & 27 & 81 \\
    SF &  & 5 & 25 & 125 & 625 \\
    AC &  & 3 & 7 & 15 & 31\\
    AC* &  & 5 & 9 & 17 & 33 \\\hline 
    \end{tabular}
\end{table}

\section{Numerical results}
\label{sec:numerics}
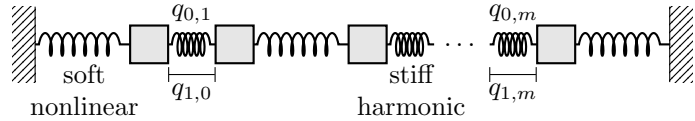
\begin{figure}[ht]
    \centering   
    \begin{tikzpicture}[every text node part/.style={align=center}]
  \def\springlength{1.25} 
  \def\springlengthcompressed{0.625}
  \tikzset{ 
      box/.style={draw,outer sep=0pt,thick}, 
      spring/.style={thick,decorate, 
        decoration={coil,pre length=0.05cm,post length=0.05cm,segment length=1.5mm, amplitude=1.5mm, aspect=0.3}}, 
      spring2/.style={thick,decorate, 
        decoration={coil,pre length=0.05cm,post length=0.05cm,segment length=0.75mm, amplitude=1.5mm, aspect=0.3}}, 
      springdescr/.style={yshift=.4cm}, 
      descriptiondescr/.style={yshift=-.6cm},
      forcedescr/.style={yshift=-.3cm},
      ground/.style={box,draw=none,fill, 
        pattern=north east lines,minimum width=0.3cm,minimum height=1cm}, 
      mass/.style={box,minimum width=.5cm,minimum height=.5cm,fill=gray!20}, 
      direction/.style={-latex,thick} 
    } 
  \node (wall) [ground] {}; 
  \draw (wall.south east) -- (wall.north east); 
 \draw [spring] (wall.east) -- node[springdescr]{} 
       +(\springlength,0)node[mass,anchor=west](M1){}; 
 \draw [spring2] (M1.east) -- node[springdescr]{$q_{0,1}$} 
        +(\springlengthcompressed,0)node[mass,anchor=west](M2){}; 
 \draw [spring] (M2.east) -- node[springdescr]{} 
        +(\springlength,0)node[mass,anchor=west](M3){}; 
 \draw [spring2] (M3.east) -- node[springdescr]{} 
        +(\springlengthcompressed,0)node[anchor=west](M4){$\dotsc$}; 
 \draw [spring2] (M4.east) -- node[springdescr]{$q_{0,m}$} 
        +(\springlengthcompressed,0)node[mass,anchor=west](Mn){}; 
\begin{scope}[xshift=0.8cm]
    \node (soft) [descriptiondescr] {soft \\ nonlinear};
\end{scope}
\begin{scope}[xshift=5.1cm]
    \node (soft) [descriptiondescr] {stiff \\ harmonic};
\end{scope}

\draw[|-|] (1.9cm,-.4cm) -- (2.525cm,-.4cm) node[draw=none,fill=none,midway,below]{$q_{1,0}$};
\draw[|-|] (6.125cm,-.4cm) -- (6.75cm,-.4cm) node[draw=none,fill=none,midway,below]{$q_{1,m}$};		
		
 \begin{scope}[xshift=8.65cm] 
 \node (wall) [ground] {}; 
 \draw (wall.south west) -- (wall.north west); 
  \draw [spring] (Mn.east) -- node[springdescr]{} 
        +(\springlength,0); 
 \end{scope} 
\end{tikzpicture}
\caption{Fermi--Pasta--Ulam (FPU) problem: $2m$ unit point masses that are connected with alternating soft nonlinear and stiff linear springs. At the end points, the springs are fixed to walls.}
\label{fig:FPU_problem}
\end{figure}

We will compare the efficiency and stability properties of the derived implicit-implicit (IMIM), implicit-explicit (IMEX) and explicit sMGARK schemes, as well as their compositions by means of the Fermi--Pasta--Ulam (FPU) problem \cite{HairerLubichWanner}, a standard benchmark problem for multirate integration of Hamiltonian systems. The FPU problem is a system that is characterized by almost-harmonic high-frequency oscillations and thus a common test problem to demonstrate the performance of structure-preserving multirate integrators (see e.g. \cite{leyendecker2013variational}). We consider a chain of $2m$ unit point masses, connected with alternating soft nonlinear and stiff linear springs, and fixed at the end points as illustrated in Fig. \ref{fig:FPU_problem}. 
 Introducing scaled displacements $q_{0,i}\ (i=1, \ldots,m)$ of the $i$th stiff spring, scaled expansion/compression $q_{1,i}$ of the $i$th stiff spring, and $p_{0,i}$, $p_{1,i}$ their velocities/momenta, the motion of the system can be described by the separable Hamiltonian system
\begin{align}\label{eq:FPU_problem}
 \begin{split}
     \H(\p,\q) = &\frac{1}{2} \sum\limits_{i=1}^m \left( p_{0,i}^2 + p_{1,i}^2\right) + \frac{\omega^2}{2} \sum\limits_{i=1}^m q_{1,i}^2 + \frac{1}{4} \bigg( (q_{0,1} - q_{1,1})^4  \\ 
     &+ \sum\limits_{i=1}^{m-1} (q_{0,i+1} - q_{1,i+1} - q_{0,i} - q_{1,i})^4 + (q_{0,m} + q_{1,m})^4 \bigg),
  \end{split}
 \end{align}
where $\q = (q_{0,1},q_{1,1},\ldots,q_{0,m},q_{1,m})^\top$, $\p = (p_{0,1},p_{1,1},\ldots,p_{0,m},p_{1,m})^\top$ and $\omega \in \mathbb{R}$ is the stiffness of the stiff springs and supposed to be large, $\omega \gg 1$. Here, the variables $q_{0,i}$, $p_{0,i}$ are slow and the variables $q_{1,i}$, $p_{1,i}$ are fast. Moreover, the third term in the Hamiltonian is the soft spring potential $\mathcal{V}^{\{\slow\}}(\q)$, while the second term is the stiff potential $\mathcal{V}^{\{\fast\}}(\q)$. The kinetic energy splits into the slow part $\mathcal{T}^{\{\slow\}}(\p) = \tfrac{1}{2}\sum_{i=1}^m p_{0,i}^2$ and the fast part $\mathcal{T}^{\{\fast\}}(\p) = \tfrac{1}{2}\sum_{i=1}^m p_{1,i}^2$. Consequently, we are able to split the Hamiltonian equations of motion into the form
\begin{align}\label{eq:FPU_EoM}
    \begin{pmatrix}\dot{\p}\\ \dot{\q}\end{pmatrix} = \begin{pmatrix}
        \mathbf{0} & -\mathbf{I} \\
        \mathbf{I} & \mathbf{0}
    \end{pmatrix} \cdot \left[ \begin{pmatrix}
        \mathcal{T}_\p^{\{\slow\}}(\p) \\ \mathcal{V}_\q^{\{\slow\}}(\q) 
    \end{pmatrix} + \begin{pmatrix}
        \mathcal{T}_\p^{\{\fast\}}(\p) \\ \mathcal{V}_\q^{\{\fast\}}(\q) 
    \end{pmatrix}\right].
\end{align}

For the simulations we consider (as in \cite{HairerLubichWanner,leyendecker2013variational}) $6$ unit point masses ($m=3$), $\omega = 50$ as the stiffness of the stiff springs, and $q_{0,1}(0) = p_{0,1}(0) = p_{1,1}(0) = 1$, $q_{1,1}(0) = \omega^{-1}$ as the initial displacements. All remaining initial values are set to zero.

\begin{figure}[H]
    \centering   
    \input{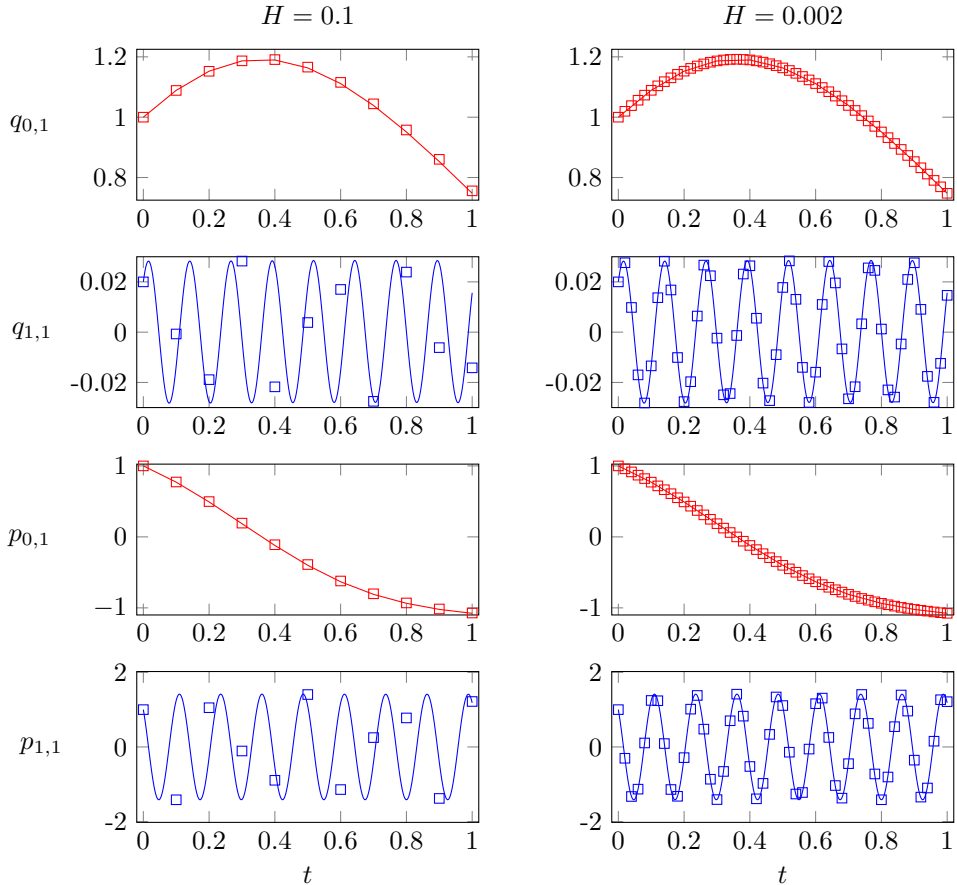}
\caption{FPU problem \eqref{eq:FPU_EoM}. Numerical approximation (denoted by the squares) for configuration and momentum of the first slow (red colored) and fast variable (blue colored) on a short time-interval $t \in [0,1]$ using the (singlerate) implicit midpoint rule. Results are shown for $H=0.1$ (left-hand side) and $H=0.002$ (right-hand side). The solid lines denote reference solutions. For visualization purposes, only every $10$th value of the numerical approximation for $H=0.002$ is displayed.}
\label{fig:FPU_singlerate_results}
\end{figure} 

\subsection{Capturing different dynamics}

\begin{figure}[ht]
    \centering
    \input{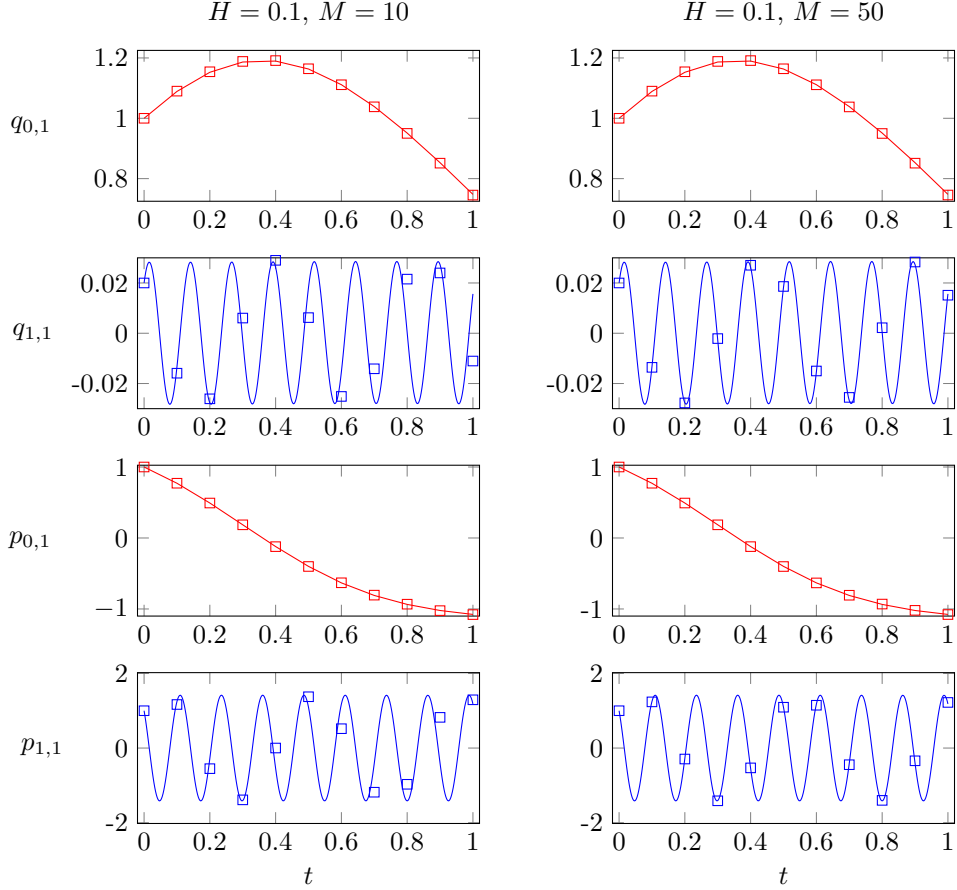}
    \caption{FPU problem \eqref{eq:FPU_EoM}. Numerical approximation (denoted by the squares) for configuration and momentum of the first slow (red colored) and fast variable (blue colored) on a short time-interval $t \in [0,1]$ using the MR-IMEX2 scheme. Results are shown for fixed macro-step size $H=0.1$ and multirate factor $M=10$ (left-hand side) and $M=50$ (right-hand side). The solid lines denote reference solutions.} 
    \label{fig:FPU_multirate_results}
\end{figure}

Figure \ref{fig:FPU_singlerate_results} shows numerical results obtained by applying the (singlerate) implicit midpoint rule. As it can be seen in the right-hand side of the figure, the step size of $H=0.002$ is small enough to resolve the fast oscillations. 
Choosing $H=0.1$ (left-hand side), the step size is small enough to capture the slow motion, however it cannot resolve the fast oscillations anymore. sMGARK methods are able to exploit this multirate behavior by keeping the macro-step size $H=0.1$ fixed while choosing a suitable multirate factor $M$ so that the micro-step size $h=H/M$ is able to capture the fast motion of the system. For $M=10$ and $M=50$, results for the MR-IMEX2 scheme \eqref{eq:IMEX2} are depicted in Figure \ref{fig:FPU_multirate_results}. The results for the MR-LPFR \eqref{eq:nested_leapfrog} and MR-IMIM2 \eqref{eq:MR-IMIM2} look similar and are not included for visualization purposes. 
The difference between the numerical approximations show that one obtains better approximations of the fast variables for increasing $M$, as expected. Furthermore, for $M=50$, the micro-step size is $h=0.002$ which is the same step size as in the right-hand side of Fig. \ref{fig:FPU_singlerate_results}. A comparison of the right-hand panels in Figure \ref{fig:FPU_singlerate_results} and  \ref{fig:FPU_multirate_results} underline that the approximation of the fast variable obtained by the multirate integrator is close to the reference solution given by the singlerate integrator, although the solution on the macro grid alone does not resolve the fast dynamics.

\subsection{Preservation of invariants}

\begin{figure}
    \centering
    \input{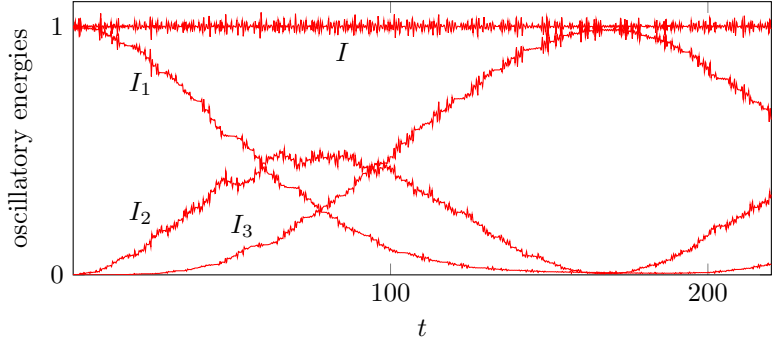}
    \caption{FPU problem \eqref{eq:FPU_EoM}. Exchange of energy in the exact solution for $t \in [0,220]$.}
    \label{fig:oscillatory_energy_exact}
\end{figure}

\begin{figure}
    \centering
    \resizebox{\linewidth}{!}{\input{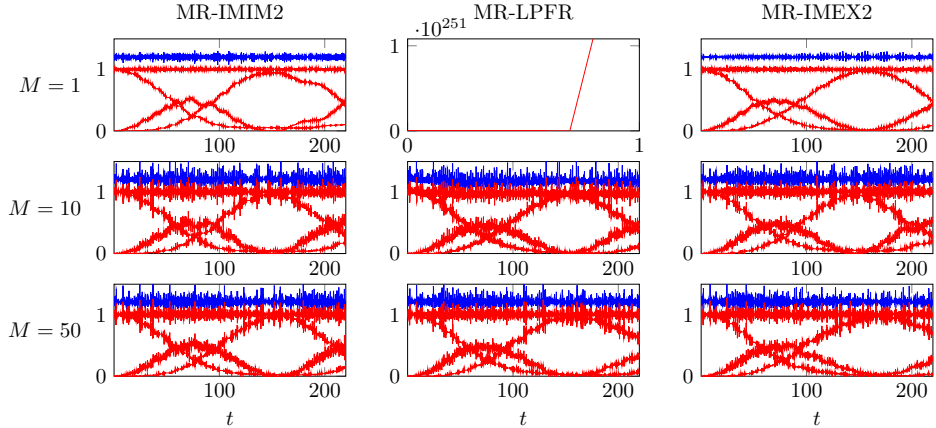}}
    \caption{FPU problem \eqref{eq:FPU_EoM}. Numerical solutions obtained on a large time interval $t \in [0,220]$ using the MR-IMIM2 (left), MR-LPFR (center) and MR-IMEX2 (right) scheme. The macro-step size is fixed at $H=0.1$; the micro-step sizes varies from row to row according to the multirate factor $M \in \{1,10,50\}$. The plots show the Hamiltonian $\H - 0.8$ (blue) and the oscillatory energies $I_1,I_2,I_3,I$ (red).}
    \label{fig:FPU_features}
\end{figure}

The exact solution of the FPU problem comes with two important features: a) the equations of motion are Hamiltonian, i.e., the total energy is exactly conserved, b) the total oscillatory energy 
\begin{align*}
    I(x,y) = \sum_{j=1}^m I_j(x_{1,j},y_{1,j}) = \frac{1}{2} \sum_{j=1}^m \left(y_{1,j}^2 + \omega^2 x_{1,j}^2\right)
\end{align*}
is an adiabatic invariant of the Hamiltonian system. More precisely, it holds 
\begin{align}
    I(x(t),y(t)) = I(x(0),y(0)) + \mathcal{O}(\omega^{-1}).
\end{align}
Fig. \ref{fig:oscillatory_energy_exact} displays the energies $I_1,I_2,I_3$ of the stiff springs and the total oscillatory energy $I=I_1+I_2+I_3$ for $t \in [0,220]$. The solution has been computed using the matlab \texttt{ode89} routine up to machine precision so that the results can be regarded as exact.
Numerical results for the oscillatory energies $I_1,I_2,I_3,I$, as well as for the Hamiltonian $\H$ are depicted in Fig. \ref{fig:FPU_features}, highlighting the good long-time behavior of the symplectic integrators. The large macro-step size $H=0.1$ exceeds the stability limit $H \omega < 2$ of the leapfrog method, explaining the stability issues in the singlerate case $M=1$. For increasing $M$, the frequency gets closer to the true behavior shown in the reference solution. However, increasing $M$ also leads to larger deviations from the exact value. Nevertheless, the values of $\H$ and $I$ stay bounded over very long time intervals.
Furthermore, the results underline the improved stability properties of the MR-IMEX2 scheme. Although it uses the leapfrog method for the slow subsystem, the critical fast subsystem is integrated using the algebraically stable midpoint-rule, resulting in an overall stable integration process that particularly leads to better results than the unconditionally stable MR-IMIM2 scheme.\\ 

\subsection{Stability of the integrators}

\begin{figure}
     \centering
        \resizebox{\linewidth}{!}
        {\input{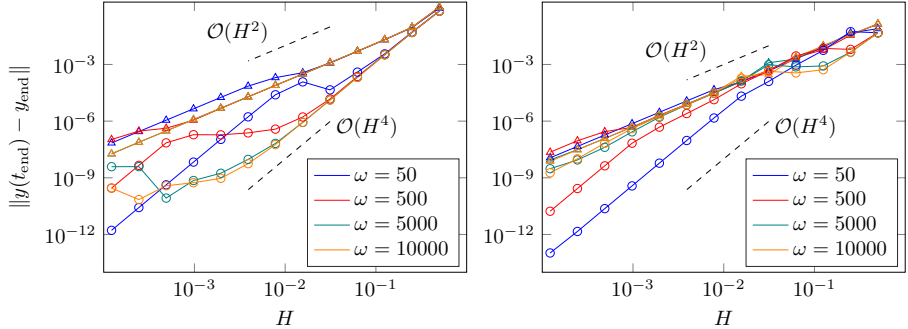}}
        \caption{FPU problem \eqref{eq:FPU_EoM}. Global error at $t_{\mathrm{end}}=3$ using the IMEX2 and IMIM2 scheme with $M=1$ (denoted by $\triangle$), as well as their compositions using the triple-jump (denoted by $\circ$). Results for varying stiffness parameter $\omega \in \{50,500,5000,10000\}$. }
        \label{fig:FPU_stability}
\end{figure}

Regarding the accuracy of the numerical simulation, we focus on the error in the slow components as the fast components will be poorly resolved in any case.  Figure \ref{fig:FPU_stability} shows the scaling of the global error at $t_{\mathrm{end}} = 3$ in the slow components with decreasing macro-step size $H$ for different values of $\omega$. The results especially point out the stability of the IMIM and IMEX scheme, as well as the stability concerns for their compositions. Both the IMIM2 and IMEX2 almost perfectly scale with the macro-step size $H$, independently of the stiffness parameter $\omega$. The use of negative step sizes in the triple-jump introduces holes in the stability domain that results in an order reduction for large values of $\omega$. Particularly, the results for the tiniest macro-step size $H=2^{-13}$ do not differ significantly for the base scheme and its triple-jump composition. As the computational cost is approximately 3 times larger for the composition scheme, the base scheme may lead to a more efficient computational process for certain $\omega$ and accuracy demands.

\section{Conclusions and outlook}
\label{sec:conclusions}
This work defines the subclass of symplectic MGARK methods as the intersection of symplectic GARK \cite{gunther2021symplectic} and MGARK methods \cite{gunther2016multirate}. The order conditions, as well as the conditions on symmetry and symplecticity are adapted to the special structure of MGARK schemes. We have shown that sMGARK schemes are algebraically stable if the weights $b_i$ are positive. Special attention is given to partitioned sMGARK schemes that are able to exploit the structure of separable Hamiltonian systems, resulting in explicit integration schemes. This introduces stability concerns since the subsystems are not stable from an analytical point of view. In this context, implicit-explicit sMGARK schemes are introduced, allowing for an explicit treatment of the slow, expensive and non-stiff part while treating the fast and stiff part implicitly, resulting in improved stability properties.
Additionally, advanced composition schemes to obtain higher-order sMGARK schemes are discussed. The required use of negative step sizes introduces holes in the stability domain so that composition schemes are not A-stable, even if the underlying sMGARK scheme is. Numerical results with the Fermi--Pasta--Ulam problem highlight the efficiency and stability properties of the schemes. 

Future work will search for symplectic and symmetric sMGARK methods of order four that do not need for negative time steps. For separable Hamiltonian systems, force-gradient integrators \cite{omelyan2002construction} are a promising tool to derive symmetric splitting methods of order four that do not need for negative time steps. Here, nested integration techniques \cite{shcherbakov2017adapted} easily extend these methods to multirate integrators. Instead of symplecticity \eqref{eq:symplecticity}, one can preserve the Hamiltonian \eqref{eq:energy_conservation} exactly. Based on averaged vector field \cite{quispel2008new} and discrete-gradient \cite{mclachlan1999geometric,gonzalez1996time} methods, we are developing energy-conserving schemes for additively partitioned Hamiltonian systems. This framework may contain operator splitting methods like \cite{frommer2023operator} as a special case and thus are of particular interest in the context of port-Hamiltonian systems.

\bmhead{Acknowledgments}
The work of K. Schäfers and M. Günther was supported by the German Research Foundation (DFG) research unit FOR5269 ``Future methods for studying confined gluons in QCD''. The work of A. Sandu was supported by the U.S. Department of Energy (DOE) 
through award ASCR DE-SC0021313, and by the Computational Science Laboratory at Virginia Tech.

\setlength{\bibsep}{0.0pt}
\def\bibfont{\scriptsize}
\bibliography{sn-bibliography}

\newpage
\begin{appendices}
\section{Order conditions for partitioned MGARK methods}
\label{sec:order-PMGARK}
\vfill
\begin{sideways}
    \small
    \rlap{\begin{minipage}{\textheight}
    \captionof{table}{Order conditions for the partitioned MGARK scheme \eqref{eq:pMGARK}-\eqref{eq:MGARK_separableHamiltonian_Butcher}}
    \label{tab:pMGARK_order}
    \begin{tabular}{l l l}
    \toprule
    $p$ & slow order condition & fast order condition \\
    \midrule
    1 & $\mathbf{\bar{b}}^{\{\slow\}\top}\mathbf{1} = 1$ & $\sum\nolimits_{\lambda=1}^M \mathbf{\bar{b}}^{\{\fast,\lambda\}\top}\mathbf{1} = M$\\
    & $\mathbf{\tilde{b}}^{\{\slow\}\top}\mathbf{1} = 1$ & $\sum\nolimits_{\lambda=1}^M \mathbf{\tilde{b}}^{\{\fast,\lambda\}\top}\mathbf{1} = M$ \\
     2 & $\mathbf{\bar{b}}^{\{\slow\}\top} \mathbf{\tilde{A}}^{\{\slow,\slow\}}\mathbf{1} = \tfrac{1}{2}$ & $\sum\nolimits_{\lambda=1}^M \mathbf{\bar{b}}^{\{\fast,\lambda\}\top} \mathbf{\tilde{A}}^{\{\fast,\fast,\lambda\}}\mathbf{1} = \frac{M}{2}$ \\
     & $\mathbf{\tilde{b}}^{\{\slow\}\top} \mathbf{\bar{A}}^{\{\slow,\slow\}}\mathbf{1} = \tfrac{1}{2}$ & $\sum\nolimits_{\lambda=1}^M \mathbf{\tilde{b}}^{\{\fast,\lambda\}\top} \mathbf{\bar{A}}^{\{\fast,\fast,\lambda\}}\mathbf{1} = \frac{M}{2}$ \\
     & $\mathbf{\bar{b}}^{\{\slow\}\top} \left(\sum_{\lambda=1}^M \mathbf{\tilde{A}}^{\{\slow,\fast,\lambda\}}\mathbf{1}\right) = \frac{M}{2}$ & $\left( \sum\nolimits_{\lambda=1}^M \mathbf{\bar{b}}^{\{\fast,\lambda\}\top} \mathbf{\tilde{A}}^{\{\fast,\slow,\lambda\}}\right)\mathbf{1} = \frac{M}{2}$ \\
     & $\mathbf{\tilde{b}}^{\{\slow\}\top} \left(\sum_{\lambda=1}^M \mathbf{\bar{A}}^{\{\slow,\fast,\lambda\}}\mathbf{1}\right) = \frac{M}{2}$ & $\left(\sum_{\lambda=1}^M \mathbf{\tilde{b}}^{\{\fast,\lambda\}\top}  \mathbf{\bar{A}}^{\{\fast,\slow,\lambda\}}\right)\mathbf{1} = \frac{M}{2}$ \\
     3 & $\mathbf{\bar{b}}^{\{\slow\}\top} \diag\!\left(\mathbf{\tilde{A}}^{\{\slow,\slow\}}\mathbf{1}\right)\mathbf{\bar{A}}^{\{\slow,\slow\}}\mathbf{1} = \frac{1}{3}$ & $\sum\nolimits_{\lambda=1}^M \mathbf{\bar{b}}^{\{\fast,\lambda\}\top} \diag\!\left(\mathbf{\tilde{A}}^{\{\fast,\fast,\lambda\}}\mathbf{1}\right)\mathbf{\bar{A}}^{\{\fast,\fast,\lambda\}}\mathbf{1} = \frac{M}{3}$ \\
     & $\mathbf{\tilde{b}}^{\{\slow\}\top} \diag\!\left(\mathbf{\bar{A}}^{\{\slow,\slow\}}\mathbf{1}\right)\mathbf{\tilde{A}}^{\{\slow,\slow\}}\mathbf{1} = \frac{1}{3}$ & $\sum\nolimits_{\lambda=1}^M \mathbf{\tilde{b}}^{\{\fast,\lambda\}\top} \diag\!\left(\mathbf{\bar{A}}^{\{\fast,\fast,\lambda\}}\mathbf{1}\right)\mathbf{\tilde{A}}^{\{\fast,\fast,\lambda\}}\mathbf{1} = \frac{M}{3}$ \\
    & $\mathbf{\bar{b}}^{\{\slow\}\top} \diag\!\left(\mathbf{\tilde{A}}^{\{\slow,\slow\}}\mathbf{1}\right) \left(\sum_{\lambda=1}^M \mathbf{\bar{A}}^{\{\slow,\fast,\lambda\}}\mathbf{1}\right) = \frac{M}{3}$ & $\left( \sum\nolimits_{\lambda=1}^M \mathbf{\bar{b}}^{\{\fast,\lambda\}\!\top\!}\!\left\{\! \diag\!\left(\mathbf{\tilde{A}}^{\{\fast,\fast,\lambda\}}\mathbf{1}\right)\! + (\lambda-1) \mathbf{I} \right\}\mathbf{\bar{A}}^{\{\fast,\slow,\lambda\}}\right) \mathbf{1} = \frac{M^2}{3}$ \\
    & $\mathbf{\tilde{b}}^{\{\slow\}\top} \diag\!\left(\mathbf{\bar{A}}^{\{\slow,\slow\}}\mathbf{1}\right) \left(\sum_{\lambda=1}^M \mathbf{\tilde{A}}^{\{\slow,\fast,\lambda\}}\mathbf{1}\right) = \frac{M}{3}$ & $\left( \sum\nolimits_{\lambda=1}^M \mathbf{\tilde{b}}^{\{\fast,\lambda\}\!\top\!}\!\left\{\! \diag\!\left(\mathbf{\bar{A}}^{\{\fast,\fast,\lambda\}}\mathbf{1}\right)\! + (\lambda-1) \mathbf{I} \right\}\mathbf{\tilde{A}}^{\{\fast,\slow,\lambda\}}\right) \mathbf{1} = \frac{M^2}{3}$ \\
     & $\mathbf{\bar{b}}^{\{\slow\}\top} \diag\!\left(\sum_{\lambda=1}^M \mathbf{\tilde{A}}^{\{\slow,\fast,\lambda\}}\mathbf{1}\right)\!\! \left(\sum_{\lambda=1}^M \mathbf{\bar{A}}^{\{\slow,\fast,\lambda\}}\mathbf{1}\right) \!=\! \frac{M^2}{3}$ & $\left(\sum_{\lambda=1}^M \mathbf{\bar{b}}^{\{\fast,\lambda\}\top} \diag\!\left( \mathbf{\tilde{A}}^{\{\fast,\slow,\lambda\}}\mathbf{1}\right)  \mathbf{\bar{A}}^{\{\fast,\slow,\lambda\}}\right)\mathbf{1} = \frac{M}{3}$ \\
     & $\mathbf{\tilde{b}}^{\{\slow\}\top} \diag\!\left(\sum_{\lambda=1}^M \mathbf{\bar{A}}^{\{\slow,\fast,\lambda\}}\mathbf{1}\right)\!\! \left(\sum_{\lambda=1}^M \mathbf{\tilde{A}}^{\{\slow,\fast,\lambda\}}\mathbf{1}\right) \!=\! \frac{M^2}{3}$ & $\left(\sum_{\lambda=1}^M \mathbf{\tilde{b}}^{\{\fast,\lambda\}\top} \diag\!\left( \mathbf{\bar{A}}^{\{\fast,\slow,\lambda\}}\mathbf{1}\right)  \mathbf{\tilde{A}}^{\{\fast,\slow,\lambda\}}\right)\mathbf{1} = \frac{M}{3}$ \\
     & $\mathbf{\bar{b}}^{\{\slow\}\top} \mathbf{\tilde{A}}^{\{\slow,\slow\}} \mathbf{\bar{A}}^{\{\slow,\slow\}} \mathbf{1} = \frac{1}{6}$ & $\sum\nolimits_{\lambda=1}^M \mathbf{\bar{b}}^{\{\fast,\lambda\}\top} \mathbf{\tilde{A}}^{\{\fast,\fast,\lambda\}} \mathbf{\bar{A}}^{\{\fast,\fast,\lambda\}} \mathbf{1} = \frac{M}{6}$ \\ 
     & $\mathbf{\tilde{b}}^{\{\slow\}\top} \mathbf{\bar{A}}^{\{\slow,\slow\}} \mathbf{\tilde{A}}^{\{\slow,\slow\}} \mathbf{1} = \frac{1}{6}$ & $\sum\nolimits_{\lambda=1}^M \mathbf{\tilde{b}}^{\{\fast,\lambda\}\top} \mathbf{\bar{A}}^{\{\fast,\fast,\lambda\}} \mathbf{\tilde{A}}^{\{\fast,\fast,\lambda\}} \mathbf{1} = \frac{M}{6}$ \\ 
     & $\mathbf{\bar{b}}^{\{\slow\}\top} \mathbf{\tilde{A}}^{\{\slow,\slow\}} \left(\sum_{\lambda=1}^M \mathbf{\bar{A}}^{\{\slow,\fast,\lambda\}}\mathbf{1}\right) = \frac{M}{6}$ & $\left( \sum\nolimits_{\lambda=1}^M \mathbf{\bar{b}}^{\{\fast,\lambda\}\top} \!\left\{\mathbf{\tilde{A}}^{\{\fast,\fast,\lambda\}} \!+\! \left( \sum\nolimits_{\mu=\lambda+1}^M \mathbf{\bar{b}}^{\{\fast,\mu\}\top} \mathbf{1}\!\right)\!\mathbf{I}\!\right\}\! \mathbf{\bar{A}}^{\{\fast,\slow,\lambda\}} \right) \!\mathbf{1} \!=\! \frac{M^2}{6}$ \\
     & $\mathbf{\tilde{b}}^{\{\slow\}\top} \mathbf{\bar{A}}^{\{\slow,\slow\}} \left(\sum_{\lambda=1}^M \mathbf{\tilde{A}}^{\{\slow,\fast,\lambda\}}\mathbf{1}\right) = \frac{M}{6}$ & $\left( \sum\nolimits_{\lambda=1}^M \mathbf{\tilde{b}}^{\{\fast,\lambda\}\top} \!\left\{\mathbf{\bar{A}}^{\{\fast,\fast,\lambda\}} \!+\! \left( \sum\nolimits_{\mu=\lambda+1}^M \mathbf{\tilde{b}}^{\{\fast,\mu\}\top} \mathbf{1}\!\right)\!\mathbf{I}\!\right\}\! \mathbf{\tilde{A}}^{\{\fast,\slow,\lambda\}} \right) \!\mathbf{1} \!=\! \frac{M^2}{6}$ \\
     & $\mathbf{\bar{b}}^{\{\slow\}\top} \left(\sum_{\lambda=1}^M  \mathbf{\tilde{A}}^{\{\slow,\fast,\lambda\}} \mathbf{\bar{A}}^{\{\fast,\slow,\lambda\}}\right)\mathbf{1} = \frac{M}{6}$ & $\sum_{\lambda=1}^M \mathbf{\bar{b}}^{\{\fast,\lambda\}\top} \mathbf{\tilde{A}}^{\{\fast,\slow,\lambda\}} \left(\sum_{\mu=1}^M \mathbf{\bar{A}}^{\{\slow,\fast,\mu\}}\mathbf{1}\right) = \frac{M^2}{6}$ \\
     & $\mathbf{\tilde{b}}^{\{\slow\}\top} \left(\sum_{\lambda=1}^M  \mathbf{\bar{A}}^{\{\slow,\fast,\lambda\}} \mathbf{\tilde{A}}^{\{\fast,\slow,\lambda\}}\right)\mathbf{1} = \frac{M}{6}$ & $\sum_{\lambda=1}^M \mathbf{\tilde{b}}^{\{\fast,\lambda\}\top} \mathbf{\bar{A}}^{\{\fast,\slow,\lambda\}} \left(\sum_{\mu=1}^M \mathbf{\tilde{A}}^{\{\slow,\fast,\mu\}}\mathbf{1}\right) = \frac{M^2}{6}$ \\
     & $\mathbf{\bar{b}}^{\{\slow\}\top} \left(\sum_{\lambda=1}^M \mathbf{\tilde{A}}^{\{\slow,\fast,\lambda\}} \{ \mathbf{\bar{A}}^{\{\fast,\fast\}} + (\lambda-1)\mathbf{I}\}\mathbf{1}\right) = \frac{M^2}{6}$ & $\left(\sum_{\lambda=1}^M \mathbf{\bar{b}}^{\{\fast,\lambda\}\top} \mathbf{\tilde{A}}^{\{\fast,\slow,\lambda\}}\right)\mathbf{\bar{A}}^{\{\slow,\slow\}}\mathbf{1} = \frac{M}{6}$ \\
     & $\mathbf{\tilde{b}}^{\{\slow\}\top} \left(\sum_{\lambda=1}^M \mathbf{\bar{A}}^{\{\slow,\fast,\lambda\}} \{ \mathbf{\tilde{A}}^{\{\fast,\fast\}} + (\lambda-1)\mathbf{I}\}\mathbf{1}\right) = \frac{M^2}{6}$ & $ \left(\sum_{\lambda=1}^M \mathbf{\tilde{b}}^{\{\fast,\lambda\}\top} \mathbf{\bar{A}}^{\{\fast,\slow,\lambda\}}\right)\mathbf{\tilde{A}}^{\{\slow,\slow\}}\mathbf{1} = \frac{M}{6}$ \\
     \botrule
    \end{tabular}%
    \vfill
    \end{minipage}}
\end{sideways}%

\newpage 
\section{Pseudocode for the multirate leapfrog scheme}
\label{sec:MR-leapfrog}

\begin{algorithm}
\caption{Multirate leapfrog \eqref{eq:nested_leapfrog}}\label{Algo:nested_Leapfrog}
\begin{algorithmic}
\State $\p \gets \p_0, \ \q \gets \q_0$
\State $\p \gets \p - \tfrac{H}{2} \mathcal{V}_\q^{\{\slow\}}\!\left(\q \right)$
\For{$\lambda=1,\ldots,M/2$}
    \State $\p \gets \p - \tfrac{h}{2} \mathcal{V}_\q^{\{\fast\}}\!\left(\q \right)$
    \State $\q \gets \q + h \mathcal{T}_\p^{\{\fast\}}\!\left(\p \right)$
    \State $\p \gets \p - \tfrac{h}{2} \mathcal{V}_\q^{\{\fast\}}\!\left(\q \right)$
\EndFor
\State $\q \gets \q + H \mathcal{T}_\p^{\{\slow\}}\!\left(\p \right)$
\For{$\lambda=1,\ldots,M/2$}
    \State $\p \gets \p - \tfrac{h}{2} \mathcal{V}_\q^{\{\fast\}}\!\left(\q \right)$
    \State $\q \gets \q + h \mathcal{T}_\p^{\{\fast\}}\!\left(\p \right)$
    \State $\p \gets \p - \tfrac{h}{2} \mathcal{V}_\q^{\{\fast\}}\!\left(\q \right)$
\EndFor
\State $\p \gets \p - \tfrac{H}{2} \mathcal{V}_\q^{\{\slow\}}\!\left(\q \right)$
\State $\p_1 \gets \p, \ \q_1 \gets \q$
\end{algorithmic}
\end{algorithm}
\end{appendices}

\end{document}